\documentclass[11pt,a4paper,final]{article}

\usepackage{amsfonts,amsmath,amssymb,theorem,graphicx, natbib}
\numberwithin{equation}{section}
\theoremstyle{plain}
\newtheorem{thm}{Theorem}[section]

\newcommand{\N}{\mathbb{N}}

\newcommand{\C}{\mathbb{C}}
\newcommand{\G}{\mathbb{G}}

\newcommand{\ind}{1\! \mathrm{l}}

\newcommand{\R}{\mathbb{R}}
\newcommand{\PP}{\mathbb{P}}

\newtheorem{cor}{Corollary}
\newtheorem{lemma}{Lemma}
\newtheorem{rem}{Remark}

\begin{document}

\title{Non parametric finite translation mixtures with dependent regime}

\author{Elisabeth  Gassiat \\
Laboratoire de Math\'ematique, Universit\'e Paris-Sud  and CNRS,
Orsay, France,
\and Judith Rousseau,\\ENSAE-CREST and CEREMADE, Universit\' e Paris-Dauphine,
Paris, France.}

\maketitle

\abstract{

In this paper we consider non parametric finite translation mixtures. We prove that all the parameters of the model are identifiable as soon as the matrix that defines the joint distribution of two consecutive latent variables is non singular and the  translation parameters are distinct. Under this assumption, we provide a consistent estimator of the number of populations, of the translation parameters and of the distribution of two consecutive latent variables, which we prove to be asymptotically normally distributed under mild dependency assumptions. We propose a non parametric estimator of the unknown translated density. In case the latent variables form a Markov chain (Hidden Markov models), we prove an oracle inequality leading to the fact that this estimator is minimax adaptive over regularity classes of densities.

}


{\bf Short title: } Non parametric finite translation mixtures \\

\section{Introduction}
Finite mixtures are widely used in applications to model heterogeneous data and to do unsupervised clustering, see for instance \citet{maclachlan:peel:2000} or \citet{marin:mengersen:robert:2004} for a review. 
Latent class models, hidden Markov models or more generally regime  switching models may be viewed as mixture models. 
Finite mixtures 
are therefore to 
be understood as  convex combinations of a finite number of probability distributions over the space the data lives in, including both static (when the latent variables are independent) and dynamical models.
Most of the developed methods use a finite dimensional description  of the probability distributions, which requires some prior knowledge of the phenomenon under investigation. 
In particular applications, it has been noticed that this may lead to poor results and various extensions have been considered. The first natural extension is to consider mixtures with an unknown number of components. This has been extensively studied  and used in the literature both from a Bayesian or frequentist point of view, see \citet{akaike:1973},  \citet{richardson:green:1997}, \citet{ishwaran:2001}, \citet{chambaz:rousseau:2008}, \citet{CGG09}, \citet{HanGas2}, to name but a few. However when the emission distribution, i.e. the distribution of each component, is misspecified this results in an overestimation of the number of components, as explained  in the discussion in \citet{rabiner:89}. Thus,  there has recently been interest in considering nonparametric mixture models in various applications, see for instance 
the discussion on the Old faithfull dataset in \citet{azzalini:bowman:90}, the need for nonparametric emission distributions in climate state identification in \citet{lambert:whiting:metcalfe:03} or the nonparametric hidden Markov model proposed in \citet{yau:papas:roberts:holmes:11}.  In absence of training data, mixture models with nonparametric emission distributions are in general not identifiable without additional structural  constraints.
 In a seminal paper, \citet{HallZhou03} discussed identifiability issues in a 2 -component nonparametric mixture model under repeated measurements (or multivarate) and showed that identifiability essentially only occured if there is at least 3 repeated measurements for each individual. This work has been extended by various authors including \citet{kasahara:shimotsu:07}, \citet{bonhomme:jochman:robin:11} and references therein. Identifiability recent results about mixtures may also be found in
\citet{AllMatRho09}.

Consider location  models 
\begin{equation}
\label{model1}
Y_{i}=m_{S_{i}}+\epsilon_{i}, \quad i \in \N
\end{equation}
where $(S_{i})_{i \in\N}$ is an unobserved sequence of random variables with finite state space $\{1,\ldots,k\}$, $(\epsilon_{i})_{i \in\N}$ is a sequence of independent identically distributed random variables taking values in $\R$,
and $m_{j}\in\R$, $j=1,\ldots,k$.
The aim is to estimate the parameters $k$, $m_{1},\ldots,m_{k}$, the distribution of the  latent variables $(S_{i})_{i\in\N}$ and the distribution $F$ of the $\epsilon_{i}$'s.
As usual for finite mixtures, one may recover the parameters only up to relabelling, and obviously, $F$ may only be estimated up to a translation (that would be reversly reported to the $m_{j}$'s).
However the identifiability issue is much more serious without further assumptions. 
 To illustrate the identifiability issues that arise with such models, assume that the $S_i $'s are independent and identically distributed. Then the $Y_i$'s are independent and have distribution
\begin{equation}
\label{model2}
P_{\mu,F} (.) = \sum_{j=1}^{k}\mu (j) F\left(\cdot - m_{j}\right).
\end{equation}
Here, $\mu(j)\geq 0$, $j=1,\ldots,k$, $\sum_{j=1}^{k}\mu\left(j\right)=1$, $m_{j}\in\R$, $j=1,\ldots,k$, and $F$ is a probability distribution on $\R$.
An equivalent representation of \eqref{model2} corresponds for instance to  $k=1$, $m_{1}=0$ and $F = P_{\mu, F}$ the marginal distribution. 
 \citet{hunter:wang:hettmanspeger:04} 
 have considered model \eqref{model2} with the additional assumption that $F$ is symmetrical and under some constraints on the $m_j$, in the case of $k \leq 4$ , see also \citet{BoMoVa:2006} and \citet{butucea:vanderkerkoeve:11} in the case where $k =2$ for an estimation procedure and asymptotic results.

In this paper, we investigate model \eqref{model1}  where the observed variables are not independent and may be non stationary. Interestingly, contrarywise to the independent case, 
we obtain identifiability without any assumption on $F$ under some very mild conditions on the process $S_1, \cdots, S_n$, see Theorem \ref{theoident}.
To be precise, if $Q$ is the $k\times k$-matrix such that $Q_{i,j}$ is the probability that $S_{1}=i$ and $S_{2}=j$, we prove that the knowledge of the distribution of $(Y_{1},Y_{2})$ allows the identification of $k$,  $m_{1},\ldots,m_{k}$, $Q$ and $F$ as soon as $Q$ is a non singular matrix, whatever $F$ may be.
Building upon our identifiability result, we propose an estimator of $k$, and of the parametric part of the distribution, namely $Q$ and $m_{1},\ldots,m_{k}$. Here, we do not need the sequence $(X_{i})_{i\in{\mathbb{N}}}$  to be strictly stationary and asymptotic stationarity is enough, then $Q$ is the stationary joint disribution of two consecutive latent variables.  
Moreover, we prove that our estimator is $\sqrt{n}$-consistent, with asymptotic Gaussian distribution, under mild dependency assumptions, see Theorem \ref{theoasymtheta}. When the number of populations is known and if the translation parameters $m_j$, $j \leq k$ are known to be bounded by a given constant, we prove that the estimator  (centered and at $\sqrt{n}$-scale) has a subgaussian distribution, see Theorem \ref{theotailstheta}.

In the context of hidden Markov models as considered in \citet{yau:papas:roberts:holmes:11}, we propose an estimator  of the non parametric part of the distribution, namely $F$, assuming that it is absolutely continuous with respect to Lebesgue measure. This estimator uses the model selection  approach developped in \citet{massart:2003}, with the penalized estimated pseudo likelihood contrast based on marginal densities
 $ \sum_{j=1}^k \hat \mu(j) f( y -\hat m_j)$.
 We prove an oracle inequality, see Theorem \ref{oracle}, which allows to deduce that our non parametric estimator is adaptive over regular classes of densities, see Theorem \ref{sadaptive} and Corollary \ref{fadaptive}.

The organization of the paper is the following. In section \ref{sec:iden} we present and prove our general identifiability theorem. In section \ref{sec:para} we define an estimator of the order and of the parametric part, and state the convergence results: asymptotic gaussian distribution and deviation inequalities. In section \ref{sec:nonpara}, we explain our non parametric estimator of the density of $F$ using model selection methods, and state an oracle inequality and adaptive convergence results. Most of the proofs are given in the Appendices.

\section{General identifiability result}
\label{sec:iden}

Let ${\cal Q}_{k}$ be the set of probability mass functions on $\{1,\ldots,k\}^{2}$, that is the set of $k\times k$ matrices $Q=(Q_{i,j})_{1\leq i,j \leq k}$ such that
for all $(i,j)\in\{1,\ldots,k\}^{2}$, $Q_{i,j}\geq 0$,
and $\sum_{i=1}^{k}\sum_{j=1}^{k}Q_{i,j}=1$.\\
We consider the joint distribution of $(Y_1, Y_2) $ under model \eqref{model1}, which has distribution 
\begin{equation}
\label{modelnota}
P_{\theta,F} (A\times B)= \sum_{i,j=1}^k Q_{ij} F(A - m_i) F(B-m_j) ,\quad  \forall A, B \in \mathcal B_{\R}
\end{equation}
where $\mathcal B_{\R}$ denotes the Borel $\sigma$ field of $\R$ and $\theta=\left(m,(Q_{i,j})_{1\leq i,j \leq k, (i,j)\neq (k,k)}\right)$, with $m=(m_{1},\ldots,m_{k})\in \R^{k}$. Recall that in this case, ordering the coefficients $m_1 \leq m_2 \leq \cdots \leq m_k$ and replacing $F$ by $F( . -m_1)$ leads to the same model so that without loss of generality we fixe $0 = m_1\leq m_2 \leq \cdots \leq m_k$.  
Let  $\Theta_{k}$ be the set of parameters $\theta$ such that 
$m_{1}=0 \leq m_{2} \leq \ldots \leq m_{k}$ and
$Q\in {\cal Q}_{k}$, where $Q=(Q_{i,j})_{1\leq i,j \leq k}$, $Q_{k,k}=1- \sum_{(i,j)\neq (k,k)}Q_{i,j}$.\\
Let also $\Theta_{k}^{0}$ be the set of parameters $\theta=\left(m,(Q_{i,j})_{1\leq i,j \leq k, (i,j)\neq (k,k)}\right) \in \Theta_{k}$ such that 
$m_{1}=0 < m_{2} < \ldots < m_{k}$ and $\rm{det}(Q)\neq 0$.
We then have the following result on the identification of $F$ and $\theta$ from $P_{\theta, F}$.
\begin{thm}
\label{theoident}
Let $F$ and $\tilde{F}$ be any probability distributions on $\R$. Let $k$ and $\tilde{k}$ be positive integers. If $\theta \in \Theta_{k}^{0}$ and $\tilde{\theta}\in \Theta_{\tilde{k}}^{0}$,
then
$$P_{\theta,F}=P_{\tilde{\theta},\tilde{F}} \Longrightarrow k=\tilde{k},\;\theta=\tilde{\theta}\;\rm{and}\;F=\tilde{F}.
$$
\end{thm}
\begin{rem}
 In the same way, it is possible to identify $\ell$-marginals, for any $\ell \geq 2$, that is the distribution of $(S_{1},\ldots,S_{\ell})$, $m$ and $F$ on the basis of the distribution of $(Y_{1},\ldots,Y_{\ell})$. 
\end{rem}

\begin{rem}
The independent case considered in \citet{hunter:wang:hettmanspeger:04},  \citet{BoMoVa:2006}, \citet{butucea:vanderkerkoeve:11} is a special case where $\det (Q) = 0$ for which our identifiability result does not hold. 
An important class of models is that of hidden Markov models. In that case, if $Q$ is the stationary distribution of two consecutive variables of the hidden Markov chain,
$\det(Q) \neq 0$ if and only if the transition matrix is non singular and the stationary distribution gives positive weights to each point. When $k=2$, we thus have $\det (Q)\neq 0$ if and only if  $S_1$ and $S_2$ are not independent.

\end{rem}

\noindent \textit{Proof of Theorem \ref{theoident} } 

Denote by $\phi_{F}$ the characteristic function of $F$, $\phi_{\tilde{F}}$ the characteristic function of $\tilde{F}$, $\phi_{\theta,1}$ (respectively $\phi_{\tilde{\theta},1}$) the characteristic function of 
the distribution of $m_{S_{1}}$ under $P_{\theta,F}$ (respectively under $P_{\tilde{\theta},\tilde{F}}$), $\phi_{\theta,2}$ (respectively $\phi_{\tilde{\theta},2}$) the characteristic function of 
the distribution of $m_{S_{2}}$ under $P_{\theta,F}$ (respectively under $P_{\tilde{\theta},\tilde{F}}$), and $\Phi_{\theta}$ (respectively $\Phi_{\tilde{\theta}}$) the characteristic function of 
the distribution of $(m_{S_{1}},m_{S_{2}})$ under $P_{\theta,F}$ (respectively under $P_{\tilde{\theta},\tilde{F}}$).
Then since the distribution of $Y_{1}$  is the same under $P_{\theta,F}$ and $P_{\tilde{\theta},\tilde{F}}$, one gets that for any $t\in \R$,
\begin{equation}
\label{chara1}
\phi_{F}\left(t\right)\phi_{\theta,1}\left(t\right)=\phi_{\tilde{F}}\left(t\right)\phi_{\tilde{\theta},1}\left(t\right).
\end{equation}
Similarly, for any $t\in \R$,
\begin{equation}
\label{chara2}
\phi_{F}\left(t\right)\phi_{\theta,2}\left(t\right)=\phi_{\tilde{F}}\left(t\right)\phi_{\tilde{\theta},2}\left(t\right).
\end{equation}
Since the distribution of $(Y_{1},Y_{2})$  is the same under $P_{\theta,F}$ and $P_{\tilde{\theta},\tilde{F}}$, one gets that for any $\mathbf t = (t_{1},t_{2})\in \R^{2}$,\begin{equation}
\label{chara3}
\phi_{F}\left(t_{1}\right)\phi_{F}\left(t_{2}\right)\Phi_{\theta}\left(\mathbf t\right)=\phi_{\tilde{F}}\left(t_{1}\right)\phi_{\tilde{F}}\left(t_{2}\right)\Phi_{\tilde{\theta}}\left(\mathbf t\right).
\end{equation}
There exists a neighborhood $V$ of $0$ such that for all $t\in V$, $\phi_{F}\left(t\right)\neq 0$, 
 so that  (\ref{chara1}), (\ref{chara2}) and (\ref{chara3}) imply that
 for any $\mathbf t  = (t_{1},t_{2})\in V^{2}$,
\begin{equation}
\label{chara5}
\Phi_{\theta}\left(\mathbf t\right)\phi_{\tilde{\theta},1}\left(t_{1}\right)\phi_{\tilde{\theta},2}\left(t_{2}\right)=\Phi_{\tilde{\theta}}\left(\mathbf t \right)\phi_{\theta,1}\left(t_{1}\right)\phi_{\theta,2}\left(t_{2}\right).
\end{equation}
Let $t_{1}$ be a fixed real number in $V$.  $\Phi_{\theta}\left(t_{1},t_{2}\right)$, $\phi_{\tilde{\theta},2}\left(t_{2}\right)$, $\Phi_{\tilde{\theta}}\left(t_{1},t_{2}\right)$, $\phi_{\theta,2}\left(t_{2}\right)$ have analytic continuations for all complex numbers $z_{2}$, $\Phi_{\theta}\left(t_{1},z_{2}\right)$, $\phi_{\tilde{\theta}}\left(z_{2}\right)$, $\Phi_{\tilde{\theta}}\left(t_{1},z_{2}\right)$, $\phi_{\theta}\left(z_{2}\right)$ which are entire functions so that (\ref{chara5}) holds with $z_{2}$ in place of $t_{2}$ for all $z_{2}$ in the complex plane $\C$ and any $t_{1}\in V$. Again, let $z_{2}$ be a fixed complex number in $\C$. $\Phi_{\theta}\left(t_{1},z_{2}\right)$, $\phi_{\tilde{\theta},1}\left(t_{1}\right)$, $\Phi_{\tilde{\theta}}\left(t_{1},z_{2}\right)$, $\phi_{\theta,1}\left(t_{1}\right)$ have analytic continuations $\Phi_{\theta}\left(z_{1},z_{2}\right)$, $\phi_{\tilde{\theta}}\left(z_{1}\right)$, $\Phi_{\tilde{\theta}}\left(z_{1},z_{2}\right)$, $\phi_{\theta}\left(z_{1}\right)$ which are entire functions so that (\ref{chara5}) holds with $z_{1}$ in place of $t_{1}$ and $z_{2}$ in place of $t_{2}$ for all $(z_{1},z_{2})\in\C^{2}$. \\
Let  now ${\cal Z}$ be the set of zeros of $\phi_{\theta,1}$, $\tilde{\cal Z}$ be the set of zeros of $\phi_{\tilde{\theta},1}$ and fix $z_{1}\in {\cal Z}$. Then, for any $z_{2}\in\C$,
\begin{equation}
\label{zeros}
\Phi_{\theta}\left(z_{1},z_{2}\right)\phi_{\tilde{\theta},1}\left(z_{1}\right)\phi_{\tilde{\theta},2}\left(z_{2}\right)=0.
\end{equation}
We now prove that $z_2 \rightarrow \Phi_{\theta}\left(z_{1},\cdot\right)$ is not the null function. For any $z\in\C$,
$$
\Phi_{\theta}\left(z_{1},z\right)=\sum_{\ell=1}^{k}\left[\sum_{j=1}^{k}Q_{\ell,j}e^{im_{j}z_{1}}
\right] e^{im_{\ell}z}.
$$
Since $0=m_{1}< m_{2}<\ldots < m_{k}$, if $\Phi_{\theta}\left(z_{1},\cdot\right)$ was the null function,  we would have for all $\ell=1,\ldots,k$
$$
\sum_{j=1}^{k}Q_{\ell,j}e^{im_{j}z_{1}}=0,
$$
which is impossible since $\det (Q) \neq 0$. Thus, $\Phi_{\theta}\left(z_{1},\cdot\right)$ is an entire function which has isolated zeros, $\phi_{\tilde{\theta},2}\left(\cdot\right)$ also, and it is possible to choose $z_{2}$ in $\C$ such that $\Phi_{\theta}\left(z_{1},z_{2}\right)\neq 0$ and
$\phi_{\tilde{\theta},2}\left(z_{2}\right)\neq 0$. Then (\ref{zeros}) leads to $\phi_{\tilde{\theta},1}\left(z_{1}\right)=0$, so that ${\cal Z} \subset \tilde{\cal Z}$.
A symmetric argument gives $\tilde{\cal Z}\subset {\cal Z}$ so that  ${\cal Z} = \tilde{\cal Z}$. Moreover, $\phi_{\theta,1}$ and $\phi_{\tilde{\theta},1}$ have growth order $1$, so that using Hadamard's factorization Theorem (see \cite{Stein:complex} Theorem 5.1) one gets that there exists a polynomial $R$ of degree $\leq 1$ such that for all $z\in\C$,
$$
\phi_{\theta,1}\left(z\right)=e^{R(z)}\phi_{\tilde{\theta},1}\left(z\right).
$$
But using $\phi_{\theta,1}\left(0\right)=\phi_{\tilde{\theta},1}\left(0\right)=1$ we get that there exists a complex number $a$ such that $\phi_{\tilde{\theta},1}\left(z\right)=e^{az}
\phi_{\theta,1}\left(z\right)$. Using now $0=m_{1}< m_{2}<\ldots < m_{k}$, and
$0=\tilde{m}_{1}< \tilde{m}_{2}<\ldots < \tilde{m}_{\tilde{k}}$ we get that $\phi_{\theta,1}=\phi_{\tilde{\theta},1}$.  Similar arguments lead to  $\phi_{\theta,2}=\phi_{\tilde{\theta},2}$.
Combining this with  (\ref{chara5}) we obtain $\Phi_{\theta}=\Phi_{\tilde{\theta}}$ which in turns implies that $k=\tilde{k}$ and $\theta=\tilde{\theta}$. Thus, using (\ref{chara1}), for all $t\in \R$ such that $\phi_{\theta,1}(t)\neq 0$,
$\phi_{F}\left(t\right)=\phi_{\tilde{F}}\left(t\right)$. Since $\phi_{\theta,1}$ has isolated zeros  and $\phi_{F}$, $\phi_{\tilde{F}}$ are continuous functions, one gets $\phi_{F}=\phi_{\tilde{F}}$ so that $F=\tilde{F}$. $\Box$

\section{Estimation of the parametric part}
\label{sec:para}

\subsection{Assumptions on the model}
Hereafter, we are given a sequence $(Y_{i})_{i\in\N}$ of real random variables with distribution $\PP^{\star}$. We assume  that (\ref{model1}) holds, with $(S_{i})_{i\in\N}$ a sequence of non-observed  random variables taking values in $\{1,\ldots,k^{\star}\}$.  We denote by $F^{\star}$ the common probability distribution of the $\epsilon_{i}$'s, and  $m^{\star}\in\R^{k^{\star}}$ the possible values of the $m_{S_{i}}$'s. We assume:
\begin{description}
\item[(A1)]
$(S_{i},S_{i+1})$ converges in distribution to $Q^{\star}\in {\cal Q}_{k^{\star}}$. \\
For $\theta^{\star}=(m^{\star},(Q^{\star}_{i,j})_{(i,j)\neq (k^{\star},k^{\star})})$, 
$\theta^{\star}\in \Theta_{k^{\star}}^{0}$,  and
all differences $m^{\star}_{j}-m^{\star}_{i}$, $i,j=1,\ldots,k^{\star}$, $i\neq j$, are distinct.
\end{description} 
We do not assume that $k^{\star}$ is known, so that the aim is to estimate $\theta^{\star}$ and $k^{\star}$ altogether. Assumption $\mathbf{(A1)}$ implies that the marginal distributions in $Q^{\star}$ are identical so that we write from now on $\phi_{\theta^\star} = \phi_{\theta^\star,1} = \phi_{\theta^\star,2}$.

The idea to estimate $\theta^{\star}$ and $k^*$ is to use equation (\ref{chara5}) which holds if and only if the parameters are equal. Consider $w$  any probability density on $\R^{2}$ with compact support $\cal S$, positive on $\cal S$ and with  $0$ belonging to the  interior of $\cal S$ ; typically $\mathcal S= [- a, a]^2$ for some positive $a$. Define, for any integer $k$ and $\theta\in\Theta_{k}$:
\begin{multline}
\label{M}
M\left(\theta\right)=\int_{\R^{2}}\left | \Phi_{\theta^{\star}}\left(t_{1},t_{2}\right) \phi_{\theta,1}\left(t_{1}\right)\phi_{\theta,2}\left(t_{2}\right)-
  \Phi_{\theta}\left(t_{1},t_{2}\right) \phi_{\theta^{\star}}\left(t_{1}\right)\phi_{\theta^{\star}}\left(t_{2}\right)   \right|^{2}\\
 \left | \phi_{F^{\star}}\left(t_{1}\right)\phi_{F^{\star}}\left(t_{2}\right) \right|^{2}
w\left(t_{1},t_{2}\right) dt_{1}dt_{2}.
\end{multline}
We shall use $M(\theta)$ as a contrast function. Indeed, thanks to Theorem \ref{theoident},  $\theta\in \Theta_{k}^{0}$ is such that $M(\theta)=0$ if and only if
$k=k^{\star}$ and $\theta=\theta^{\star}$. \\
We estimate  $M(\cdot)$ by
\begin{equation}
\label{Mn}
M_{n}\left(\theta\right)=\int_{\R^{2}}\left | \widehat{\Phi}_{n}\left(t_{1},t_{2}\right) \phi_{\theta,1}\left(t_{1}\right)\phi_{\theta,2}\left(t_{2}\right)-
  \Phi_{\theta}\left(t_{1},t_{2}\right) \widehat{\phi}_{n,1}\left(t_{1}\right)\widehat{\phi}_{n,2} \left(t_{2}\right) \right|^{2} 
w\left(t_{1},t_{2}\right) dt_{1}dt_{2},
\end{equation}
where $\widehat{\Phi}_{n}$ is an estimator of the characteristic function of the asymptotic distribution of $(Y_{t},Y_{t+1})$, 
$\widehat{\phi}_{n,1}(t)=\widehat{\Phi}_{n}(t,0)$ and $\widehat{\phi}_{n,2}(t)=\widehat{\Phi}_{n}(0,t)$. One may take for instance the empirical estimator
\begin{equation}
\label{estiphi}
\widehat{\Phi}_{n}\left(t_{1},t_{2}\right)=\frac{1}{n}\sum_{j=1}^{n-1} \exp i\left(t_{1}Y_{j}+t_{2}Y_{j+1}\right).
\end{equation}
We require that $\widehat{\Phi}_{n}$ is uniformly upper bounded;
if  $\widehat{\Phi}_{n}$  is defined by \eqref{estiphi} then  it is uniformly upper bounded by $1$.
Define, for any $\mathbf t =(t_{1},t_{2}) \in\R^{2}$
$$
Z_{n}\left(\mathbf t\right)=
\sqrt{n}\left( \widehat{\Phi}_{n}(\mathbf t)  - \Phi_{\theta^{\star}}(\mathbf t)\phi_{F^{\star}}\left(t_{1}\right)\phi_{F^{\star}}\left(t_{2}\right)\right).
$$
Our main assumptions on the model and on the estimator  $\widehat{\Phi}_{n}$ are the following. 
\begin{description}
\item[(A2)]
The process $(Z_{n}\left(\mathbf t\right))_{\mathbf t \in {\cal S}}$ converges weakly to a Gaussian process $(Z\left(\mathbf t\right))_{\mathbf t \in {\cal S}}$
 in the set of complex continuous functions on $\cal S$ endowed with the uniform norm and with covariance kernel $\Gamma(\cdot,\cdot)$.
\end{description}

\begin{description}
\item[(A3)] There exist real numbers $E$ and $c$ (depending on $\theta^{\star}$)  such that for all $x\geq 0$ and $n\geq 1$,
$$
\PP^{\star}\left(\sup_{\mathbf t \in {\cal S}}\left| Z_{n}\left(\mathbf t\right)\right| \geq E + x \right) \leq \exp \left(-cx^{2} \right).
$$
\end{description}

$(\mathbf{A2})$ will be used to obtain the asymptotic distribution of the estimator, and $(\mathbf{A3})$ to obtain non asymptotic deviation inequalities.
Note that  $(\mathbf{A2})$ and $(\mathbf{A3})$ are for instance verified if we use (\ref{estiphi}),  under stationarity and mixing conditions on the $Y_{j}$'s. This follows applying results of \citet{doukhan:massart:rio:1994}, \citet{doukhan:massart:rio:1995} and \citet{rio:2000}.

\subsection{Definition of the estimator} \label{subsec:estpara}
Our contrast function verifies $M\left(\theta\right)=0$ if and only if $\theta=\theta^{\star}$ only when we restrict  $\theta$ to belong to $\cup_{k\in\N}\Theta_{k}^{0}$.
When minimization is performed over $\cup_{k\in\N}\Theta_{k}^{0}$ it may happen that the minimizer is on the boundary. To get rid of this problem, we build our estimator $\widehat{\theta}_{n}$ using a preliminary consistent estimator $\tilde{\theta}_{n}$, and then restrict the minimization using the information given by $\tilde{\theta}_{n}$.\\
Define for any integer $k$, $I_{k}$ a positive continuous function on $\Theta_{k}^{0}$ and tending to $+\infty$ on the boundary of $\Theta_{k}^{0}$ or whenever $\|m\|$ tends to infinity. For instance one may take
$$
I_{k}\left(m,(Q_{i,j})_{(i,j)\neq(k,k)}\right)= -\log \rm{det} Q - \sum_{i=2}^{k}\log \frac{|m_{i}-m_{i-1}|}{(1+\|m\|_{\infty})^{2}}.
$$
Let $(k_{n},\tilde{\theta}_{n})$ be a minimizer over $\{(k,\theta):k\in\N, \theta\in\Theta_{k}\}$ of
$$
C_{n}\left(k,\theta\right)=M_{n}\left(\theta\right) + \lambda_{n}\left[J\left(k\right) + I_{k}\left(\theta\right)\right]
$$
where $J:\N\rightarrow \N$ is an increasing function tending to infinity at infinity and $(\lambda_{n})_{n\in\N}$ a decreasing sequence of real numbers tending to $0$ at infinity such that
\begin{equation}
\label{lambdan}
\lim_{n\rightarrow+\infty}\sqrt{n}\lambda_{n}=+\infty
\end{equation}

Define now $\widehat{\theta}_{n}$ as a minimizer of $M_{n}$ over 
$$
\left\{ \theta\in \Theta_{k_n}\;:\;  I_{k_n}\left(\theta\right)\leq 2 I_{k_n}\left(\tilde{\theta}_{n}\right)\right\}.
$$ 

In case $k^{\star}$ is known, we may choose another estimator. Let ${\cal K}$ be a compact subset of $\Theta_{k^{\star}}^{0}$. We denote by $\overline{\theta}_{n}(\cal K)$ a minimizer of $M_{n}$ over $\cal K$. This estimator will also be used 
as a theoretical trick in the proof of the asymptotic distribution of $\widehat{\theta}_{n}$. 

\subsection{Asymptotic results}
Our first  result gives the asymptotic distribution of $\widehat{\theta}_{n}$.
To define the asymptotic variance, we define $\nabla M \left(\theta\right)$ the gradient of $M$ at point $\theta$ and $D_{2}M\left(\theta\right)$ the Hessian of $M$ at point $\theta$. We also set 
$V$  the variance of the gaussian process
\begin{multline*}
\int \left\{C\left(\mathbf t\right)
 \left[ Z\left(-\mathbf t\right)\phi_{\theta^{\star}}\left(-t_{1}\right)\phi_{\theta^{\star}}\left(-t_{2}\right)-\Phi_{\theta^{\star}}\left(-\mathbf t\right)\left(Z(-t_{1},0)\phi_{\theta^{\star}}\left(-t_{2}\right)+Z(0,-t_{2})\phi_{\theta^{\star}}\left(-t_{1}\right)\right) \right]\right.\\
\left.+ C\left(-\mathbf t\right)
\left[ Z \left(\mathbf t\right)\phi_{\theta^{\star}}\left(t_{1}\right)\phi_{\theta^{\star}}\left(t_{2}\right)-\Phi_{\theta^{\star}}\left(\mathbf t\right)\left(Z(t_{1},0)\phi_{\theta^{\star}}\left(t_{2}\right)+Z(0,t_{2})\phi_{\theta^{\star}}\left(t_{1}\right)\right) \right]
\right\}w\left(\mathbf t\right)d\mathbf t
\end{multline*}
where
$$
C\left(\mathbf t\right)= \Phi_{\theta^{\star}}\left(\mathbf t\right)\nabla\left(  \phi_{\theta^{\star}}\left(t_{1}\right)\phi_{\theta^{\star}}\left(t_{2}\right)\right)-\nabla\Phi_{\theta^{\star}}\left(\mathbf t\right)\phi_{\theta^{\star}}\left(t_{1}\right)\phi_{\theta^{\star}}\left(t_{2}\right).
$$

\begin{thm}
\label{theoasymtheta}
Assume $(\mathbf{A1})$, $(\mathbf{A2})$, and (\ref{lambdan}). Then $D_{2}M\left(\theta^{\star}\right)$ is non singular, and
for any compact subset ${\cal K}$ of $\Theta_{k^{\star}}^{0}$ such that $\theta^{\star}$ lies in the interior of ${\cal K}$,  $\sqrt{n} (\overline{\theta}_{n}(\cal K)-\theta^*)$ converges in distribution to the centered Gaussian with variance
$$
\Sigma =D_{2}M\left(\theta^{\star}\right)^{-1} V D_{2}M\left(\theta^{\star}\right)^{-1}.
$$
Moreover,  $\sqrt{n} (\widehat{\theta}_{n}-\theta^*)$ converges in distribution to the centered Gaussian with variance 
 $\Sigma$.
\end{thm}

If one wants to use Theorem \ref{theoasymtheta} to build confidence sets, one needs to have a consistent estimator of $\Sigma$. Since $D_{2}M$ is a continuous functions of $\theta$, $D_{2}M\left(\widehat{\theta}_{n}\right)$ is a consistent estimator of $D_{2}M\left(\theta^{\star}\right)$.  Also, $V$ may be viewed as a continuous function of $\Gamma (\cdot,\cdot)$ and $\theta$, as easy but tedious computations show. One may use empirical estimators of $\Gamma (\cdot,\cdot)$  which are uniformly consistent under stationarity and mixing conditions, to get a consistent estimator of $V$. This leads to a plug-in consistent estimator of $\Sigma$. \\
Another possible way to estimate $\Sigma$ is  to use a boostrap method, following for instance \citet{clemencon:garivier:tressou:09} when the hidden variables form a Markov chain.

When we have deviation inequalities for the process $Z_{n}$, we are able to provide deviation inequalities for $\sqrt{n} (\overline{\theta}_{n}(\cal K)-\theta^*)$. Such inequalities have interest by themselves, they will also be used for proving adaptivity of our non parametric estimator in Section \ref{sec:nonpara}. 

\begin{thm}
\label{theotailstheta}
Assume $(\mathbf{A1})$ and $(\mathbf{A3})$. Let ${\cal K}$ be a compact subset of $\Theta_{k^{\star}}^{0}$ such that $\theta^{\star}$ lies in the interior of ${\cal K}$.  
Then there exist real numbers $c^{\star}$,  $M^{\star}$,
and an integer $n^{\star}$ such that for all $n\geq n^{\star}$ and $M\geq M^{\star}$,
$$
\PP^{\star}\left( \sqrt{n}\|\overline{\theta}_{n}({\cal K})-\theta^{\star}\|\geq M \right)\leq 8\exp\left(-c^{\star}M^{2}\right).
$$
In particular, for any integer $p$,
$$
\sup_{n\geq 1}E_{\PP^{\star}} \left[\left(\sqrt{n}\|\overline{\theta}_{n}({\cal K})-\theta^{\star}\| \right)^{p}\right] < +\infty.$$
\end{thm}

\section{Estimation of the non parametric part in the case of  hidden Markov models} 
\label{sec:nonpara}

In this section we assume that $\PP^{\star}$ is the distribution of a stationary ergodic hidden Markov model (HMM for short), that is the sequence
$(S_{t})_{t\in{\mathbb{N}}}$ is a stationary ergodic Markov chain. We also assume that 
 the unknown distribution $F^{\star}$ has density $f^{\star}$ with respect to Lebesgue measure. Thus the density $s^{\star}$ of $Y_{1}$ writes 
 $$
 s^{\star}\left(y\right)=\sum_{j=1}^{k^{\star}}\mu^{\star}\left(j\right)f^{\star}\left(y-m_{j}^{\star}\right),
 $$
 where $\mu^{\star}(j)=\sum_{i=1}^{k^{\star}}Q^{\star}_{j,i}$, $1\leq i \leq k^{\star}$. We shall assume moreover:
 \begin{description}
\item[(A4)]
For all $i,j=1,\ldots,k^{\star}$, $Q^{\star}_{i,j}>0$, and there exists $\delta >0$ such that
$$
\int_{\R} \left[f^{\star}\left(y\right)\right]^{1-\delta} dy < +\infty.
$$
\end{description}
Notice that, if the observations form a stationary HMM  and if for all $i,j=1,\ldots,k^{\star}$, $Q^{\star}_{i,j}>0$, then the sequence is geometrically uniformly ergodic, and applying results of  \citet{doukhan:massart:rio:1994}, \citet{doukhan:massart:rio:1995} and \citet{rio:2000}, $\mathbf{(A2)}$ and  $\mathbf{(A3)}$ hold if we use \eqref{estiphi}.\\

We propose to use model selection methods to estimate $f^{\star}$ using penalized marginal likelihood. 
We assume in this section that $k^{\star}$ is known, and that we are given an estimator $\widehat{\theta}_{n}=((\widehat{m}_{i})_{1\leq i \leq k^{\star}}, (\widehat{Q}_{i,j})_{(i,j)\neq (k^{\star},k^{\star})})=\overline{\theta}_{n}({\cal K})$ of $\theta^{\star}$ for some compact subset ${\cal K}$ of $\Theta_{k^{\star}}^{0}$ such that $\theta^{\star}$ lies in the interior of ${\cal K}$.
Let $\widehat{\mu}(i)=\sum_{j=1}^{k^{\star}}\widehat{Q}_{i,j}$,
{$1\leq i \leq k^{\star}$. 
Define for any density function $f$ on $\R$ 
$$
\ell_{n}\left(f\right)=\frac{1}{n}\sum_{i=1}^{n}\log \left[ \sum_{j=1}^{k^{\star}} \widehat{\mu}\left(j\right)f\left(Y_{i}-\widehat{m}_{j}\right) \right].
$$
Let ${\cal F}$ be the set of probability densities on $\R$. We shall use the model collection $({\cal F}_{p})_{p\geq 2}$ 
of Gaussian mixtures with $p$ components as approximation of ${\cal F}$.
Let us define for any integer $p$
\begin{equation}
\label{sieve}
{\cal F}_{p}=\left\{\sum_{i=1}^{p}\pi_{i} \varphi_{u_{i}}\left(x-\alpha_{i}\right), \;\alpha_{i}\in [-A_{p},A_{p}],\;u_{i}\in [b_{p}, B],\; \pi_{i}\geq 0,\;i=1,\ldots,p,\;
\sum_{i=1}^{p}\pi_{i}=1
\right\}
\end{equation}
where $B$ and $A_{p}, b_{p}$, $p\geq 2$, are positive real numbers, and where $\varphi_{\beta}$ is the Gaussian density with variance $\beta^{2}$ given by $\varphi_{\beta}(x)=\exp (-x^{2}/2\beta^{2}) / \beta \sqrt{2\pi}$.
For any $p\geq 2$, let $\widehat{f}_{p}$ be the maximizer of $\ell_{n}(f)$ over ${\cal F}_{p}$. Define
$$
D_{n}\left(p\right)=-\ell_{n}\left(\widehat{f}_{p} \right)+\text{pen}\left( p,n \right).
$$
Our model selection estimator $\widehat{f}$ will be given by $\widehat{f}_{\widehat{p}}$ whenever $\widehat{p}$ is a minimizer of $D_{n}$. 

\subsection{Oracle inequality}
The following theorem says that a suitable choice of the penalty term $\text{pen}\left( p ,n\right)$ leads to an estimator having good non asymptotic and asymptotic properties.
In the following, 
$$
\widehat{s}_{\widehat{p}}\left(\cdot \right)=\sum_{j=1}^{k^{\star}} \widehat{\mu}\left(j\right)\widehat{f}_{\widehat{p}}\left(\cdot-\widehat{m}_{j}\right),
$$
is the estimator of $s^{\star}$,
$$
S_{p}^{\star}=\{\sum_{j=1}^{k^{\star}}\mu^{\star}\left(j\right)f\left(\cdot-m^{\star}_{j}\right), f\in {\cal F}_{p}  \}
$$ 
for any $p\geq 2$, $h^{2}(\cdot,\cdot)$ is the Hellinger distance and $K(\cdot,\cdot)$  the Kullback-Leibler divergence between probability densities.
For any $p\geq 1$, fix some $f_{p}\in {\cal F}_{p}$ and set $s_{p}=\sum_{j=1}^{k^{\star}}\mu^{\star}\left(j\right)f_{p}\left(\cdot-m^{\star}_{j}\right)$. Of course to derive good behaviour of the estimator from the oracle inequality, one will have to  choose carefully $f_{p}$. 
\begin{thm}
\label{oracle}
Assume $\mathbf{(A1)} $, $\mathbf{(A3)}$ and $\mathbf{(A4)}$. Let $(x_{p})_{p\geq 2}$ be a sequence of positive real numbers such that $\Sigma=\sum_{p\geq 2}e^{-x_{p}} < +\infty$.
Then there exist positive real numbers $\kappa$ and $C$, depending only on $Q^{\star}$ and $\delta$ such that, as soon as
$$
\text{pen}\left(p,n\right) \geq \frac{\kappa}{n} \left( k^{\star}p \left[\log n + \log\left(\frac{1}{b_{p}}\right)+\log A_{p}\right]+x_{p}\left[1+\left|\log \left(1+\frac{1}{b_{p}^{\delta}}\right)\right|\right]\right),
$$
one has
$$
E_{\PP^{\star}}\left[h^{2}\left(s^{\star},\widehat{s}_{\widehat{p}}\right)\right] \leq 
C \left\{\inf_{p\geq 2}\left(K\left(f^{\star},f_{p}\right)+\text{pen}\left(p,n\right)+ E_{\PP^{\star}}\left[V_{p}\right] 
\right)+\frac{\Sigma}{n} 
\right\}
$$
with
$$
V_{p}=
\frac{1}{n}\sum_{i=1}^{n} \log \left(\frac{\sum_{j=1}^{k^{\star}}\widehat{\mu}(j)f_{p}(Y_{i}-\widehat{m}_{j})}{\sum_{j=1}^{k^{\star}}\mu^{\star}(j)f_{p}(Y_{i}-m^{\star}_{j})} \right).
$$
\end{thm}
The proof of Theorem \ref{oracle} is postponed to Appendix \ref{appb}.\\
Notice that the constant in the so-called oracle inequality depends on $\PP^{\star}$, so that the result of Theorem \ref{oracle} is not of real practical use. 
Also, the upper bound depends on  $\widehat{\theta}$, for which the results in Section \ref{sec:para} are for large enough $n$.
However,  Theorem \ref{oracle} is the 
building stone to understand how to choose a penalty function and to prove adaptivity of our estimator. 

\subsection{Adaptive estimation}

We prove now that $\widehat{s}_{\widehat{p}}$ is an adaptive estimator of $s^{\star}$, and that, if $\max_{j}\mu^{\star}(j)>\frac{1}{2}$,  $\widehat{f}_{\widehat{p}}$ is an adaptive estimator of $f^{\star}$. Adaptivity will be proved on  the following classes of regular densities.

Let $y_{0}>0$, $c>0$, $M>0$, $\tau >0$, $C>0$, $\lambda>0$ and $L$ a positive polynomial function on $\R$. Let also $\beta >0$ and $\gamma > (3/2-\beta)_+$.
If we denote ${\mathcal P}=(y_{0},c_{0},M,\tau,C,\lambda,L)$, we define  $\mathcal H_{loc}(\beta, \gamma, {\cal P})$ as the set of probability densities $f$ on $\R$ satisfying:
\begin{itemize}
\item $f$ is monotone on $(-\infty , -y_0)$ and on $(y_0, +\infty )$, and $\inf_{|y| \leq y_0 }f(y) \geq c_0>0$.
\item 
 \begin{equation} \label{cond:tail}
  \forall y\in\R,\; f(y) \leq M e^{- \tau |y|}
 \end{equation}
 \item
 $\log f$ is $\lfloor \beta \rfloor$  times continuously differentiable with derivatives $\ell_j$, $j\leq \beta$ satisfying for all $x\in \R$ and all $|y-x|\leq \lambda$,
 \begin{equation*} \label{loc-H}
  |\ell_{\lfloor\beta\rfloor}(y)  - \ell_{\lfloor\beta\rfloor}(x) | \leq \lfloor\beta\rfloor ! L(x) |y-x|^{\beta - \lfloor\beta\rfloor}
  \end{equation*}
 and 
 \begin{equation*}\label{moment}
 \int_{\R} |\ell_j(y)|^{\frac{2 \beta + \gamma }{ j } }f (y) dy \leq C.
 \end{equation*}
\end{itemize}
We use $\widehat{s}_{\widehat{p}}$ where the penalty is set to
 $$
\text{pen}\left(p,n\right) = \frac{3\kappa}{n}( k^{\star}p + x_{p})\log n.
$$
\begin{thm}
\label{sadaptive}
Assume $\mathbf{(A1)} $, $\mathbf{(A3)} $ and $\mathbf{(A4)}$. Then for any $\cal P$, $\beta \geq 1/2$ and $\gamma > (3/2-\beta)_+$, there exists $C(\beta, \gamma, {\cal P})>0$ such that
$$
\limsup_{n\rightarrow +\infty} \left(\frac{n}{ (\log n)^3}\right)^{\frac{2\beta}{2\beta +1}}\sup_{f^{\star}\in\mathcal H_{loc}(\beta, \gamma, {\cal P}) }E_{\PP^{\star}}\left[h^{2}\left(s^{\star},\widehat{s}_{\widehat{p}}\right)\right] \leq 
C(\beta, \gamma, {\cal P}).
$$
\end{thm}

Thus, $\widehat{s}_{\widehat{p}}$ is adaptive on the regularity $\beta$ of the density classes up to $(\log n)^{3\beta/(2\beta + 1)}$, see \citet{maugis:michel} for a lower bound of the asymptotic minimax risk in the case of independent and identically distributed random variables. Using Theorem \ref{sadaptive}, we can  also derive adaptive asymptotic rates for the
minimax $L_1$-risk  for the estimation of $f^*$.

\begin{cor}
\label{fadaptive}
Assume $\mathbf{(A1)} $, $\mathbf{(A3)}$, $\mathbf{(A4)}$ and that  $\max_{j}\mu^{\star}(j)>\frac{1}{2}$ . Then for any $\cal P$, $\beta \geq 1/2$ and $\gamma > (3/2-\beta)_+$, 
\begin{equation*}
\limsup_{n\rightarrow +\infty} \left(\frac{n}{ (\log n)^3}\right)^{\frac{\beta}{2\beta +1}}\sup_{f^{\star}\in\mathcal H_{loc}(\beta, \gamma, {\cal P}) }E_{\PP^{\star}}\left[\left\|\widehat{f}_{\widehat{p}}-f^{\star}\right\|_{1}\right] 
\leq 
\frac{2\sqrt{C(\beta, \gamma, {\cal P})} }{\left(2 \max_{j}\mu^{\star}\left(j\right)-1\right)}.
\end{equation*}
\end{cor}

It is possible that the constraint, $\max_j \mu^{\star}(j) >1/2$ is not sharp, however note that the Fourier transform of $s^{\star}$ is expressed as 
 $\phi_{\theta^{\star}} \phi_{f^{\star}}$ with 
  $ \phi_{\theta^{\star}}(t) = \sum_{j=1}^k \mu^{\star}(j) e^{it m^{\star}_j}$ and $\phi_{f^{\star}}$ the Fourier transform of $f^{\star}$, and that $|\phi_{\theta^{\star}}(t)|>0$ for all $t\in \R$ if and only if $\max_j \mu^{\star}(j) >1/2$, applying the main theorem of \citet{moreno:73}. 

\noindent \textit{ Proof of Corollary \ref{fadaptive}}

We shall use 
$$
\|s^{\star}-\widehat{s}_{\widehat{p}}\|_{1} \leq 2 h\left(s^{\star},\widehat{s}_{\widehat{p}}\right),
$$
together with
\begin{equation*}
\begin{split}
\|s^{\star}-\widehat{s}_{\widehat{p}}\|_{1}&= \|\sum_{j=1}^{k^{\star}}\mu^{\star}\left(j\right)f^{\star}\left(\cdot-m_{j}^{\star}\right)-
\sum_{j=1}^{k^{\star}} \widehat{\mu}\left(j\right)\widehat{f}_{\widehat{p}}\left(\cdot-\widehat{m}_{j}\right)\|_{1}\\
&\geq \| \sum_{j=1}^{k^{\star}} \mu^{\star}\left(j\right)(\widehat{f}_{\widehat{p}}-f^{\star})\left(\cdot-\widehat{m}_{j}\right)\|_{1}
 -\|\widehat{\theta}_{n}-\theta^{\star}\|\\
& \quad -\|\sum_{j=1}^{k^{\star}}\mu^{\star}\left(j\right)\left(f^{\star}\left(\cdot-m_{j}^{\star}\right)
-f^{\star}\left(\cdot-\widehat{m}_{j}\right)\right)\|_{1}
\\
&\geq \left(2 \max_{j}\mu^{\star}\left(j\right)-1\right)\left\|\widehat{f}_{\widehat{p}}-f^{\star}\right\|_{1}
-\|\widehat{\theta}_{n}-\theta^{\star}\|\\
& \quad -\|f^{\star}\left(\cdot-m_{j}^{\star}\right)
-f^{\star}\left(\cdot-\widehat{m}_{j}\right)\|_{1}
\end{split}
\end{equation*}
which follows by using iteratively the triangle inequality. Using $\beta\geq 1/2$, Theorem \ref{theotailstheta}.
and Theorem \ref{sadaptive}, we thus get that
\begin{equation*}
\limsup_{n\rightarrow +\infty} \left(\frac{n}{ (\log n)^3}\right)^{\frac{\beta}{2\beta +1}}\sup_{f^{\star}\in\mathcal H_{loc}(\beta, \gamma, {\cal P}) }E_{\PP^{\star}}\left[\left\|\widehat{f}_{\widehat{p}}-f^{\star}\right\|_{1}\right] 
\leq 
\frac{2\sqrt{C(\beta, \gamma, {\cal P})} }{\left(2 \max_{j}\mu^{\star}\left(j\right)-1\right)}
\end{equation*}
as soon as
\begin{equation}
\label{lemmecor}
\lim_{n\rightarrow +\infty} \left(\frac{n}{ (\log n)^3}\right)^{\frac{\beta}{2\beta +1}}\sup_{f^{\star}\in\mathcal H_{loc}(\beta, \gamma, {\cal P}) }
E_{\PP^{\star}}\left[\|f^{\star}\left(\cdot-m_{j}^{\star}\right)
-f^{\star}\left(\cdot-\widehat{m}_{j}\right)\|_{1}\right]=0.
\end{equation}
Now, since $f^{\star}\in H_{loc}(\beta, \gamma, {\cal P})$ with $\beta \geq 1/2$, if  $|\widehat m_j - m_j^{\star}|\leq \lambda$,
$$|\log f^{\star}(y - \widehat m_j)- \log f^{\star}(y - m_j^{\star}) | \leq L(y-m_j^{\star}) |\widehat m_j - m_j^{\star}|^{\beta \wedge 1}.$$
Set $M\geq \frac{1}{2c^{\star}}$, and $a>0$ such that, if $|y|\leq n^{a}$, then $L(y) |\widehat m_j - m_j^{\star}|^{\beta \wedge 1}\leq 1$.
Observe also that since $\widehat{\theta}_{n}$ stays in a compact set, for large enough $n$, if $|y|\geq n^{a}$, then for any $j$,
$|y-\widehat m_j|\geq n^{a}/2$ and $|y-m_j^{\star}|\geq n^{a}/2$.
We obtain, using $|e^{u}-1|\leq 2u$ for $0\leq u \leq 1$:
\begin{eqnarray*}
\|f^{\star}\left(\cdot-m_{j}^{\star}\right)
-f^{\star}\left(\cdot-\widehat{m}_{j}\right)\|_{1}& \leq &
2 \left(\frac{M\log n}{n}\right)^{-(\beta \wedge 1)/2} \int L(y-m_j^{\star})f^{\star}(y-m_{j}^{\star})dy\\
&& + 2 \int_{|y|\geq n^{a}/2} f^{\star}(y)dy +  \ind_{\|\theta^{\star}-\widehat \theta_n\| > \sqrt{M\log n}/\sqrt{n}},
\end{eqnarray*}
and \eqref{lemmecor} follows from  Theorem \ref{theotailstheta}, $\beta \geq 1/2$ and the  fact that $f^{\star}\in H_{loc}(\beta, \gamma, {\cal P})$ has exponentially decreasing tails.
$\Box$

\subsection{Computation of $\widehat{f}_{p}$}
The computation of $\widehat{f}_{p}$ may be performed using the EM-algorithm, which is particularly simple for Gaussian mixtures. Indeed, for
$f=\sum_{i=1}^{p}\pi_{i} \varphi_{\beta_{i}}\left(\cdot-\alpha_{i}\right)$, $\sum_{j=1}^{k^{\star}} \widehat{\mu}\left(j\right)f\left(\cdot-\widehat{m}_{j}\right)$ is a mixture
of $pk^{\star}$ Gaussian densities $\phi_{\beta_{i}}\left(\cdot-\alpha_{i}-\widehat{m}_{j}\right)$ with weights $\pi_{i} \widehat{\mu}\left(j\right)$. Starting from an initial point
$((\pi_{i}^{0})_{1\leq i \leq p},(v_{i}^{0})_{1\leq i \leq p},(\sigma_{i}^{0})_{1\leq i \leq p})$, the EM $l$-th iteration may be easily computed as
$$
\pi_{i}^{l+1}=\frac{\sum_{j=1}^{k^{\star}}\sum_{t=1}^{n} \widehat{\mu}\left(j\right)\pi_{i}^{l}\varphi_{\beta_{i}^{l}}\left(Y_{t}-\widehat{m}_{j}-\alpha_{i}^{l}\right)}{\sum_{i'=1}^{p}\sum_{j=1}^{k^{\star}}\sum_{t=1}^{n} \widehat{\mu}\left(j\right)\pi_{i'}^{l}\varphi_{\beta_{i'}^{l}}\left(Y_{t}-\widehat{m}_{j}-\alpha_{i'}^{l}\right)},\;i=1,\ldots,p,
$$
$$
\alpha_{i}^{l+1}=T_{-A_{p},A_{p}}\left[\frac{\sum_{j=1}^{k^{\star}}\sum_{t=1}^{n}\left(Y_{t}-\widehat{m}_{j}\right) \widehat{\mu}\left(j\right)\pi_{i}^{l}\varphi_{\beta_{i}^{l}}\left(Y_{t}-\widehat{m}_{j}-\alpha_{i}^{l}\right)}{\sum_{j=1}^{k^{\star}}\sum_{t=1}^{n}\widehat{\mu}\left(j\right)\pi_{i}^{l}\varphi_{\beta_{i}^{l}}\left(Y_{t}-\widehat{m}_{j}-\alpha_{i}^{l}\right)} \right]
,\;i=1,\ldots,p,
$$
where for any real numbers $C_{1},C_{2}$, $T_{C_{1},C_{2}}$ is the troncature function: $T_{C_{1},C_{2}}(x)=x\ind_{C_{1}\leq x \leq C_{2}}+C_{1}\ind_{x<C_{1}}+C_{2}\ind_{x>C_{2}}$, and
$$
\sigma_{i}^{l+1}=T_{b_{p}, B}\left[\frac{\sum_{j=1}^{k^{\star}}\sum_{t=1}^{n}\left(Y_{t}-\widehat{m}_{j}-v_{i}^{l}\right)^{2} \widehat{\mu}\left(j\right)\pi_{i}^{l}\varphi_{\beta_{i}^{l}}\left(Y_{t}-\widehat{m}_{j}-\alpha_{i}^{l}\right)}{\sum_{j=1}^{k^{\star}}\sum_{t=1}^{n}\widehat{\mu}\left(j\right)\pi_{i}^{l}\varphi_{\beta_{i}^{l}}\left(Y_{t}-\widehat{m}_{j}-\alpha_{i}^{l}\right)} \right]
,\;i=1,\ldots,p.
$$

\section*{Acknowledgements}
This work was partly supported by  the 2010--2014 grant ANR Banhdits AAP Blanc SIMI 1.

\appendix

\section{Proof of Theorem 3.1}
First of all, we prove a lemma we shall use several times.
Using $| |A|^{2}-|B|^{2}| \leq |A-B| ||A|+|B||$ and the fact that characteristic functions are uniformly upper bounded by $1$, we get that for any integer $k$ and any $\theta\in\Theta_{k}$:
\begin{multline*}
\left| M_{n}\left(\theta\right)- M\left(\theta\right)\right| \leq 2 \int \left\{\left| \widehat{\Phi}_{n}\left(t_{1},t_{2}\right) -\Phi_{\theta^{\star}}\left(t_{1},t_{2}\right)\phi_{F^{\star}}\left(t_{1}\right)\phi_{F^{\star}}\left(t_{2}\right)\right|\right.\\
\left.+
\left|\widehat{\phi}_{n}\left(t_{1}\right)\widehat{\phi}_{n} \left(t_{2}\right)- \phi_{\theta^{\star}}\left(t_{1}\right)\phi_{\theta^{\star}}\left(t_{2}\right)\phi_{F^{\star}}\left(t_{1}\right)\phi_{F^{\star}}\left(t_{2}\right)
\right|\right\}w\left(t_{1},t_{2}\right) dt_{1}dt_{2}.
\end{multline*}
The upper bound does not depend on $k$ and $\theta$, $ \widehat{\Phi}_{n}$ is uniformly upper bounded, and we get
\begin{equation}
\label{majo}
\sup_{k\geq 2,\;\theta\in\Theta_{k}}\left| M_{n}\left(\theta\right)- M\left(\theta\right)\right| =O\left( \sup_{\mathbf t \in {\cal S}}\left|\frac{ Z_{n}(\mathbf t)}{\sqrt{n}}\right|\right)=O_{\PP^{\star}}(1/\sqrt{n})
\end{equation}
which together with Theorem \ref{theoident} gives
\begin{lemma}
\label{proptool}
If $(k_{n},\theta_{n})_{n}$, $\theta_{n}\in\Theta_{k_{n}}$,  is a random sequence such that there exists an integer $K\geq k^*$, and a compact subset ${\cal T}$ of
$\cup_{k\leq K}\Theta_{k}^{0}$ such that
$$\PP^*\left(k_{n}\leq K\;\rm{and}\;\theta_{n}\in{\cal T}\right) \rightarrow 1\;\rm{and}\;M_{n}\left(\theta_{n}\right)=o_{\PP^{\star}}(1),
$$
then
$$\PP^*\left(k_{n}=k^{\star}\right) \rightarrow 1\;\rm{and}\;\theta_{n}=\theta^{\star} + o_{\PP^{\star}}(1).
$$
\end{lemma}
Since $C_{n}\left(k_{n},\tilde{\theta}_{n}\right)\leq C_{n}\left(k^*,\theta^*\right)$ and $M_{n}$ is a non negative function, we get
$$
\left[J(k_{n}) + I_{k_{n}}\left(\tilde{\theta}_{n}\right)\right] \leq \left[J(k^{\star}) + I_{k^{\star}}\left(\theta^{\star}\right)\right]+\frac{M_{n}\left(\theta^{\star}\right)-M\left(\theta^{\star}\right)}{\lambda_{n}},
$$
so that using (\ref{majo}), assumption $\mathbf{(A2)}$ and (\ref{lambdan}) we get
\begin{equation}
\label{negli}
\left[J(k_{n} )+ I_{k_{n}}\left(\tilde{\theta}_{n}\right)\right] \leq \left[J(k^{\star}) + I_{k^{\star}}\left(\theta^{\star}\right)\right]+o_{\PP^{\star}}\left(1\right).
\end{equation}
Also,
$$
M_{n}\left(\tilde{\theta}_{n}\right)\leq M_{n}\left(\theta^{\star}\right) + \lambda_{n}\left[J(k^{\star}) + I_{k^{\star}}\left(\theta^{\star}\right)\right],
$$
so that
$$
M_{n}\left(\tilde{\theta}_{n}\right)=o_{\PP^{\star}}\left(1\right).
$$
Thus, using  (\ref{negli}) and Lemma \ref{proptool}
\begin{equation}
\label{resuinter1}
\PP^*\left(k_{n}=k^{\star}\right) \rightarrow 1\;\rm{and}\;\tilde{\theta}_{n}=\theta^{\star} + o_{\PP^{\star}}(1).
\end{equation}
Set now 
${\cal K}=\left\{ \theta\in \Theta_{k^{\star}}\;:\;  I_{k^{\star}}\left(\theta\right)\leq 4 I_{k^{\star}}\left(\theta^{\star}\right)\right\}$. $\cal K$ is a compact subset of $\Theta_{k^{\star}}^{0}$. 
Let  $E_{n}$ be the event $(k_{n}=k^{\star}\;{\rm{and}}\;\widehat{\theta}_{n}=\overline{\theta}_{n}({\cal K}))$. Using Lemma  \ref{proptool}, we get that $\overline{\theta}_{n}(\cal K)
$ is a consistent estimator of $\theta^{\star}$, and using (\ref{resuinter1}) and Lemma  \ref{proptool}, we get also that $\widehat{\theta}_{n}
$ is a consistent estimator of $\theta^{\star}$, so  $M_n$ has the same minimizer on $\mathcal K$ and on $\{ I_{k_n}(\theta) \leq 2 I_{k_n}( \tilde \theta_n) \}$, with probability tending to $1$, since it belongs to a neigbourhood of $\theta^*$. Thus, 
$
\PP^{\star}\left( E_{n}\right) \rightarrow 1.
$
Now, since
$$
\widehat{\theta}_{n}=\overline{\theta}_{n}({\cal K})\ind_{E_n} + \widehat{\theta}_{n}\ind_{E_n^c},
$$
Theorem \ref{theoasymtheta} follows as soon as we prove that  $\sqrt{n} (\overline{\theta}_{n}({\cal K})-\theta^{\star})$ converges in distribution to the centered Gaussian with variance $\Sigma$. But this is a straighforward consequence of 
$$
D_{2}M_{n}\left(\theta_{n}\right)\left(\overline{\theta}_{n}({\cal K})-\theta^{\star}\right)=\nabla M_{n} \left(\theta^{\star}\right),
$$
for some  $\theta_{n}\in\Theta_{k^{\star}}$ such that $\|\theta_{n}-\theta^{\star}\| \leq \|\overline{\theta}_{n}({\cal K})-\theta^{\star}\| $, the consistency of $\overline{\theta}_{n}({\cal K})$  and the following Lemma
\begin{lemma}
\label{gradetd2}
Assume {\bf (A1)} and {\bf (A2)}. Then
\begin{itemize}
\item
$\sqrt{n}\nabla M_{n} \left(\theta^{\star}\right)$ converges in distribution to a centered gaussian with variance $V$.
\item $D_{2}M\left(\theta^{\star}\right)$ is non singular, 
and
for any random variable $\theta_{n}\in\Theta_{k^{\star}}$ converging in $\PP^{\star}$-probability to $\theta^{\star}$, one has
$$
D_{2}M_{n}\left(\theta_{n}\right)=D_{2}M\left(\theta^{\star}\right)+o_{\PP^{\star}}\left(1\right).
$$
\end{itemize}
\end{lemma}

\noindent \textit{Proof of Lemma \ref{gradetd2}}

First notice that, in every formula, taking the conjugate of any involved function at point $\mathbf t$ is the same as taking the function  at point $-\mathbf t$. This is also verified for derivatives.
Write now for any $\theta\in\Theta_{k^{\star}}$ and any $\mathbf t =(t_{1},t_{2})$
$$
G_{n}\left(\theta,\mathbf t\right)= \widehat{\Phi}_{n}\left(\mathbf t\right) \phi_{\theta,1}\left(t_{1}\right)\phi_{\theta,2}\left(t_{2}\right)-
  \Phi_{\theta}\left(\mathbf t\right) \widehat{\phi}_{n,1}\left(t_{1}\right)\widehat{\phi}_{n,2} \left(t_{2}\right)
$$
so that, if $\nabla G_{n}\left(\theta,\mathbf t\right)$ denotes the gradient of $G_{n}$ with respect to $\theta$ at point $(\theta,\mathbf t)$, one has
$$
\nabla M_{n} \left(\theta^{\star}\right)=\int \left[\nabla G_{n}\left(\theta^{\star},\mathbf t\right)G_{n}\left(\theta^{\star},-\mathbf t\right)+
\nabla G_{n}\left(\theta^{\star},-\mathbf t\right)G_{n}\left(\theta^{\star},\mathbf t\right)\right] w\left(\mathbf t\right)d\mathbf t.
$$
Now, writing $ \widehat{\Phi}_{n}\left(\mathbf t\right)=\frac{Z_{n}(\mathbf t)}{\sqrt{n}}+ \Phi_{\theta^{\star}}(\mathbf t)\phi_{F^{\star}}(t_{1})\phi_{F^{\star}}(t_{2})$ and using ${\bf (A2)}$ one gets easily
\begin{multline*}
\sqrt{n}\nabla M_{n} \left(\theta^{\star}\right)=\int \left\{\phi_{F^{\star}}(t_{1})\phi_{F^{\star}}(t_{2})\left[ \Phi_{\theta^{\star}}\left(\mathbf t\right)\nabla\left(  \phi_{\theta^{\star}}\left(t_{1}\right)\phi_{\theta^{\star}}\left(t_{2}\right)\right)-\nabla\Phi_{\theta^{\star}}\left(\mathbf t\right)\phi_{\theta^{\star}}\left(t_{1}\right)\phi_{\theta^{\star}}\left(t_{2}\right)\right]\right.\\
\left. \left[ Z_{n}\left(-\mathbf t\right)\phi_{\theta^{\star}}\left(-t_{1}\right)\phi_{\theta^{\star}}\left(-t_{2}\right)-\Phi_{\theta^{\star}}\left(-\mathbf t\right)\left(Z_{n}(-t_{1},0)\phi_{\theta^{\star}}\left(-t_{2}\right)+Z_{n}(0,-t_{2})\phi_{\theta^{\star}}\left(-t_{1}\right)\right) \right]\right.\\
\left.+\phi_{F^{\star}}(-t_{1})\phi_{F^{\star}}(-t_{2})\left[ \Phi_{\theta^{\star}}\left(-\mathbf t\right)\nabla\left(  \phi_{\theta^{\star}}\left(-t_{1}\right)\phi_{\theta^{\star}}\left(-t_{2}\right)\right)-\nabla\Phi_{\theta^{\star}}\left(-\mathbf t\right)\phi_{\theta^{\star}}\left(-t_{1}\right)\phi_{\theta^{\star}}\left(-t_{2}\right)\right]\right.\\
\left. \left[ Z_{n}\left(\mathbf t\right)\phi_{\theta^{\star}}\left(t_{1}\right)\phi_{\theta^{\star}}\left(t_{2}\right)-\Phi_{\theta^{\star}}\left(\mathbf t\right)\left(Z_{n}(t_{1},0)\phi_{\theta^{\star}}\left(t_{2}\right)+Z_{n}(0,t_{2})\phi_{\theta^{\star}}\left(t_{1}\right)\right) \right]
\right\}w\left(\mathbf t\right)d\mathbf t\\
+O_{\PP^{\star}}\left(\frac{1}{\sqrt{n}}\right)
\end{multline*}
and the convergence in distribution of $\sqrt{n}\nabla M_{n} \left(\theta^{\star}\right)$ to a centered gaussian with variance $V$ follows.

Similar computation gives that for any  $\theta\in\Theta_{k^{\star}}$
  \begin{multline*}
  D_{2} M_n(\theta) -D_{2} M_n(\theta^{\star})=
  \int | \widehat \Phi_n(\mathbf t) |^2 \left[ A_{1}(\mathbf t,\theta)-A_{1}(\mathbf t,\theta^{\star})\right]
  w(\mathbf t)d\mathbf t \\
 + \int | \widehat \Phi_n(t_1,0) |^2| \widehat \Phi_n(0,t_2) |^2 \left[ A_{2}(\mathbf t,\theta)-A_{2}(\mathbf t,\theta^{\star})\right]
  w(\mathbf t)d\mathbf t \\
 +  Re\left\{\int  \widehat \Phi_n(-\mathbf t)  \widehat \Phi_n(t_1,0) \widehat \Phi_n(0,t_2)
   \left[ A_{3}(\mathbf t,\theta)-A_{3}(\mathbf t,\theta^{\star})\right]
   w(\mathbf t)d\mathbf t   \right\}
  \end{multline*}
for matrix-valued functions $A_{1}(\mathbf t,\theta)$, $A_{2}(\mathbf t,\theta)$, $A_{3}(\mathbf t,\theta)$ that are, in a neighborhood of $\theta^{\star}$, continuous in the variable $\theta$ for all $\mathbf t$ and uniformly upper bounded. Thus $D_{2} M_n({\theta}_n) -D_{2} M_n(\theta^{\star})$ converges in $\PP^{\star}$-probability to $0$ whenever ${\theta}_n$ is a random variable converging in  $\PP^{\star}$-probability to $\theta^{\star}$. 

Finally, note that at point $\theta^{\star}$ the Hessian of $M$ simplifies into:
\begin{equation*}
D_{2} M(\theta^{\star}) = 2 \int   H(\mathbf t)   H(-\mathbf t)^{T} \left | \phi_{F^{\star}}\left(t_{1}\right)\phi_{F^{\star}}\left(t_{2}\right) \right|^{2}w(\mathbf t) d\mathbf t,
\end{equation*}
with 
$$H(\mathbf t) =  \Phi_{\theta^{\star}}( \mathbf t)(\phi_{\theta^{\star}}(t_1) \nabla \phi_{\theta^{\star}}(t_2) + \nabla \phi_{\theta^{\star}}(t_1) \phi_{\theta^{\star}}(t_2)) - \nabla \Phi_{\theta^{\star}}(\mathbf t) \phi_{\theta^{\star}}(t_1)\phi_{\theta^{\star}}(t_2).$$
Denote by $H_{m_j}(\mathbf t)$, $j=2,\ldots,k^{\star}$, $H_{Q_{j_1,j_2}}(\mathbf t)$, $j_{1},j_{2}=1,\ldots,k^{\star}$, $(j_1,j_2)\neq (k^{\star},k^{\star})$ the components of the vector
$H(\mathbf t)$.
Definite positiveness of the second derivative of $M$ at $\theta^{\star}$ can thus be established by proving that, if for all ${\mathbf{t}}\in {\cal S}$, 
 \begin{equation} \label{ind:U}
 \sum_{j=2}^k U_{m_j} H_{m_j}(\mathbf t) + \sum_{(j_1,j_2)\neq (k,k)} U_{j_1,j_2}H_{Q_{j_1,j_2}}(\mathbf t) = 0
 \end{equation}
then
  $$U_{m_j}=0, \;j=2,\cdots,k^{\star}, \;U_{j_1,j_2}=0, \;j_{1},j_{2}=1,\ldots,k^{\star},\;(j_1,j_2)\neq (k^{\star},k^{\star}).$$
By linear independence of the functions $e^{it a}$ and $t e^{itb}$ this implies in particular that for all $\mathbf t =(t_{1},t_{2})$,
\begin{multline} \label{ind:mj}
\sum_{j_1,\cdots, j_4=1}^{k^{\star}} U_{m_{j_1}}\mu^{\star}\left(j_1\right)\mu^{\star}\left(j_2\right)Q^{\star}_{j_3,j_4}e^{it_1(m^{\star}_{j_1}+m^{\star}_{j_3})+it_2(m^{\star}_{j_2}+m^{\star}_{j_4})} \\
=\sum_{j_1,\cdots, j_4=1}^{k^{\star}} U_{m_{j_1}}\mu^{\star}\left(j_2\right)\mu^{\star}\left(j_3\right)Q^{\star}_{j_1,j_4}e^{it_1(m^{\star}_{j_1}+m^{\star}_{j_3})+it_2(m^{\star}_{j_2}+m^{\star}_{j_4})}
\end{multline}
with $U_{m_1}=0$. The smallest possible term $m^{\star}_{j_1}+m^{\star}_{j_3}$ with $j_1>1$ is equal to $m^{\star}_2 = m^{\star}_2+m^{\star} _1$ setting $j_1=2$ and $j_3=1$ only. Thus \eqref{ind:mj} implies that 
$$U_{m_2}\mu^{\star}\left(2\right)\sum_{j_2,j_4=1}^{k^{\star}} \mu^{\star}\left(j_2\right)Q^{\star}_{1,j_4}e^{it_2(m^{\star}_{j_2}+m^{\star}_{j_4})} = U_{m_2} \mu^{\star}\left(1\right)\sum_{j_2, j_4=1}^{k ^{\star}}\mu^{\star}\left(j_2\right)Q^{\star}_{2,j_4}e^{it_2(m^{\star}_{j_2}+m^{\star}_{j_4})}$$
for all $t_2$, i.e.
$$U_{m_2}\mu^{\star}\left(2\right)\phi_{\theta^{\star}}(t_2)\sum_{j_4=1}^{k^{\star}} Q^{\star}_{1,j_4}e^{it_2m^{\star}_{j_4}} = U_{m_2} \mu^{\star}\left(1\right)\phi_{\theta^{\star}}(t_2)\sum_{j_4=1}^{k^{\star}} Q^{\star}_{2,j_4}e^{it_2m^{\star}_{j_4}}.$$
Since $\phi_{\theta^{\star}} $ has only isolated  zeros this is satisfied if and only if 
$$U_{m_2}\mu^{\star}\left(2\right)\sum_{j_4=1}^{k^{\star}} Q^{\star}_{1,j_4}e^{it_2m^{\star}_{j_4}} = U_{m_2} \mu^{\star}\left(1\right)\sum_{j_4=1}^{k^{\star}} Q^{\star}_{2,j_4}e^{it_2m^{\star}_{j_4}}.$$
Thus \eqref{ind:mj} is satisfied only if either $U_{m_2}=0$ or
$\mu^{\star}\left(2\right)Q^{\star}_{1,j}=\mu^{\star}\left(1\right)Q^{\star}_{2,j}$ for all $j$.
The latter is impossible since $Q^{\star}$ is non singular, thus $U_{m_2}=0$ and \eqref{ind:mj} becomes
\begin{multline*}
\sum_{j_1=3,j_{2},\cdots, j_4=1}^{k^{\star}} U_{m_{j_1}}\mu^{\star}\left(j_1\right)\mu^{\star}\left(j_2\right)Q^{\star}_{j_3,j_4}e^{it_1(m^{\star}_{j_1}+m^{\star}_{j_3})+it_2(m^{\star}_{j_2}+m^{\star}_{j_4})} \\
=\sum_{j_1=3,j_{2},\cdots, j_4=1}^{k^{\star}} U_{m_{j_1}}\mu^{\star}\left(j_2\right)\mu^{\star}\left(j_3\right)Q^{\star}_{j_1,j_4}e^{it_1(m^{\star}_{j_1}+m^{\star}_{j_3})+it_2(m^{\star}_{j_2}+m^{\star}_{j_4})}
\end{multline*}
The smallest possible value for $m^{\star}_{j_1}+m^{\star}_{j_3}$ is then $m^{\star}_3$ which is obtained with the only configuration $j_1=3, j_3=1$. The same argument as before leads to $U_{m_3}=0$. Iteration of the argument leads to
$U_{m_j}= 0$ for all $j=1, \cdots , k^{\star}$. We now study the derivatives associated to $Q$.
We write $U$ the $k^{\star}\times k^{\star}$-matrix whose components are $U_{j_1,j_2}$ for $(j_1,j_2) \neq (k^{\star},k^{\star})$ and $U_{k^{\star},k^{\star}} = -\sum_{(j_1,j_2)\neq (k^{\star},k^{\star})} U_{j_1,j_2}$. Then
$$\sum_{(j_1,j_2)\neq (k^{\star},k^{\star})} U_{j_1,j_2} \nabla_{Q_{j_1,j_2}} \Phi_{\theta^{\star}}(\mathbf t) = V(t_1)^{T} UV(t_2)$$
where for any $t\in\R$, $V(t) = ((e^{itm^{\star}_{j}})_{j=1, \cdots, k^{\star}})^{T}$,
and
$$\sum_{(j_1,j_2)\neq (k^{\star},k^{\star})} U_{j_1,j_2} \nabla_{Q_{j_1,j_2}} \phi_{\theta^{\star}}( t_1) = V(t_1)^{T} U \ind
$$
with $\ind = (1,\cdots , 1)^{T}\in \R^{k^{\star}}$, since  $\phi_{\theta^{\star}}(t_1) = V(t_1)^{T}Q^{\star}\ind $ and $\Phi_{\theta^{\star}}( \mathbf t) = V(t_1)^T Q^{\star} V(t_2)$.
We can then express \eqref{ind:U} as
\begin{multline}
\label{equaQ}
V(t_1)^{T} \left[Q^{\star} V(t_2) V(t_2)^{T} U\ind \ind ^{T}(Q^{\star})^{T}+ Q^{\star} V(t_2) V(t_2)^{T} Q^{\star}\ind \ind^{T}U^{T}\right.\\
\left.
 - UV(t_2) V(t_2)^{T}Q\ind \ind ^{T}(Q^{\star})^{T}
 \right]V(t_1)= 0.
\end{multline}
 Note also that since all differences $m^{\star}_{j_{1}}-m^{\star}_{j_{2}}$, $j_{1}\neq j_{2}$, are distinct, if $A$ is a $k^{\star}\times k^{\star}$-matrix and $\cal I$ is an open subset of $\R$,
\begin{equation} 
\label{truc}
\left[\forall t \in {\cal I},\;V(t)^{T}AV(t)=0\right] \Longrightarrow A+A^{T}=0.
\end{equation}
Then \eqref{equaQ} implies
 \begin{multline} \label{equality:V2}
 Q^{\star}V(t_2)V(t_2)^{T} U\ind \ind^{T} (Q^{\star})^{T}+Q^{\star}\ind \ind^{T}U^{T} V(t_2)V(t_2)^{T} (Q^{\star})^{T}\\
 +Q^{\star} V(t_2) V(t_2)^{T} Q^{\star}\ind \ind^{T} U^{T}
 + U\ind \ind^{T}(Q^{\star})^{T}  V(t_2) V(t_2)^{T} (Q^{\star})^{T}\\ 
  - UV(t_2)V(t_2)^{T}Q^{\star}\ind \ind^{T}(Q^{\star})^{T}
   - Q^{\star}\ind \ind^{T}(Q^{\star})^{T}V(t_2)V(t_2)^{T}U^{T}=0.
 \end{multline} 
 Recall also that $\ind^{T} U \ind = 0$ and that $Q^{\star}\ind = \mu^{\star}$.  Note that $U\ind = \alpha \mu^{\star} $ with $\alpha \in \R$ if and only if $\alpha=0$ since $\ind^{T} U \ind = 0$ while $\ind^{T} \mu^{\star}=1$. Therefore if $U\ind \neq 0$ there exists $w\in \R^{k^{\star}}$ such that $w^{T} (U\ind) \neq 0$ while $(\mu^{\star})^{T} w=0$. Multiplying the above equality on the left by $w^{T} $ and on the right by $w$ leads to
$$ w^{T} Q^{\star} V(t_2) V(t_2)^{T}(\mu^{\star}) (U\ind)^{T} w =0
$$
that for all $t_{2}$ in an open set. Using \eqref{truc} again and since  $(U\ind)^{T} w\neq 0$ we get that 
$$ \mu^{\star} [(Q^{\star})^{T} w]^{T}  + [(Q^{\star})^{T} w]( \mu^{\star})^{T} = 0. $$
Since $\mu^{\star}(j) >0$ for all $j$ this implies that 
 $ (Q^{\star})^{T} w= 0$ which is impossible since $Q^{\star}$ has full rank. Therefore $U\ind= 0$ and \eqref{equality:V2} becomes
 $V(t_2)^{T}\mu^{\star} [UV(t_2)( \mu^{\star})^{T} + \mu^\star V(t_2)^T U^T ] = 0$,  
 that is
 $UV(t_2)(\mu^{\star})^{T} +\mu^\star V(t_2)^T U^T=0$ for all $t_{2}$ in an open set. Multiplying on the left by $\ind$ implies that
 $UV(t_2) =0$ for all $t_2$ in an open set so that $U=0$. 
 $\Box$
 
 
 \section{Proof of Theorem 3.2}

Define for any $\theta\in\Theta_{k^{\star}}$,
$L_{n}(\theta)=M_{n}(\theta)-M(\theta)$. Then, since $M_{n}(\overline{\theta}_{n}({\cal K}))\leq M_{n}(\theta^{\star})$, one easily gets
$$
M\left(\overline{\theta}_{n}({\cal K})\right)-M\left(\theta^{\star}\right) \leq \left|L_{n}\left(\overline{\theta}_{n}({\cal K})\right)-L_{n}\left(\theta^{\star}\right)\right|.
$$
Define for any $\mathbf t =(t_{1},t_{2})$ and any $\theta$
$$
G\left(\theta,\mathbf t\right)=\left\{ \Phi_{\theta^{\star}}\left(\mathbf t\right) \phi_{\theta,1}\left(t_{1}\right)\phi_{\theta,2}\left(t_{2}\right)-
  \Phi_{\theta}\left(\mathbf t\right) \phi_{\theta^{\star},1}\left(t_{1}\right)\phi_{\theta^{\star},2} \left(t_{2}\right)\right\} \phi_{F^{\star}}\left(t_{1}\right)
   \phi_{F^{\star}}\left(t_{2}\right)
$$
and
\begin{multline*}
B_{n}\left(\theta,\mathbf t\right)=\phi_{F^{\star}}\left(t_{1}\right)
   \phi_{F^{\star}}\left(t_{2}\right)\left\{\frac{Z_{n}(\mathbf t)}{\sqrt{n}}\phi_{\theta,1}\left(t_{1}\right)\phi_{\theta,2}\left(t_{2}\right)\right.\\
\left.-\Phi_{\theta}\left(\mathbf t\right)
\left[ \frac{Z_{n}(t_{1},0)}{\sqrt{n}}\phi_{\theta,2}\left(t_{2}\right)+\frac{Z_{n}(0,t_{2})}{\sqrt{n}}\phi_{\theta,1}\left(t_{1}\right)+  \frac{Z_{n}(t_{1},0)Z_{n}(0,t_{2})}{n} \right]
\right\} 
\end{multline*}
Writing $ \widehat{\Phi}_{n}\left(\mathbf t\right)=\frac{Z_{n}(\mathbf t)}{\sqrt{n}}+ \Phi_{\theta^{\star}}(\mathbf t)\phi_{F^{\star}}(t_{1})\phi_{F^{\star}}(t_{2})$ one gets
$$
L_{n}\left(\theta\right)=\int \left(\left[ B_{n}\left(\theta,\mathbf t\right)+G\left(\theta,\mathbf t\right)\right]
\left[ B_{n}\left(\theta,-\mathbf t\right)+G\left(\theta,-\mathbf t\right) \right]  - |G\left(\theta,\mathbf t\right)|^2\right) w\left(\mathbf t\right)d\mathbf t.
$$
Since $G\left(\theta^{\star},\mathbf t\right)=0$ for all $\mathbf t$ we obtain
\begin{multline*}
L_{n}\left(\theta\right)-L_{n}\left(\theta^{\star}\right)=\int \left\{|B_{n}\left(\theta,\mathbf t\right)|^{2}-|B_{n}\left(\theta^{\star},\mathbf t\right)|^{2}+
B_{n}\left(\theta,\mathbf t\right)G\left(\theta,-\mathbf t\right)\right.\\
\left.+B_{n}\left(\theta,-\mathbf t\right)G\left(\theta,\mathbf t\right)\right\}w\left(\mathbf t\right)d\mathbf t
\end{multline*}
which gives
\begin{multline*}
\left|L_{n}\left(\theta\right)-L_{n}\left(\theta^{\star}\right)\right|\leq \int \left\{\left|B_{n}\left(\theta,\mathbf t\right)-B_{n}\left(\theta^{\star},\mathbf t\right)\right|\left|B_{n}\left(\theta,\mathbf t\right)+B_{n}\left(\theta^{\star},\mathbf t\right)\right|\right.\\
\left. +2 \left|B_{n}\left(\theta,\mathbf t\right)\right| \left|G\left(\theta,\mathbf t\right)-G\left(\theta^{\star},\mathbf t\right)\right|\right\}w\left(\mathbf t\right)d\mathbf t
\end{multline*}
which leads to
\begin{equation}
\label{majo1}
M\left(\overline{\theta}_{n}({\cal K})\right)-M\left(\theta^{\star}\right)\leq C W_{n}\|\overline{\theta}_{n}({\cal K})-\theta^{\star}\|
\end{equation}
for some constant $C$ and any integer $n$, and with 
$$W_{n}=\left\{\frac{V_{n}}{\sqrt{n}}+\frac{V_{n}^{2}}{n}+\frac{V_{n}^{3}}{n^{3/2}}+\frac{V_{n}^{4}}{n^{2}}\right\},
\;V_{n}=\sup_{\mathbf t \in {\cal S}}\left| Z_{n}\left(\mathbf t\right)\right|.
$$
 Observe now that, since $D_{2}M$ is continuous and $D_{2}M(\theta^{\star})$ is non singular, there exists $\lambda >0$ and $\alpha >0$ such that, if $\|\theta-\theta^{\star}\|\leq \alpha$, then $M(\theta)-M(\theta^{\star}) \geq \frac{\lambda}{2} \|\theta-\theta^{\star}\|^{2}$. Moreover, there exists $\delta >0$ such that, if $\theta\in {\cal K}$ is such that 
 $\|\theta-\theta^{\star}\| \geq \alpha$, then $M(\theta)-M(\theta^{\star})\geq \delta$. Using \eqref{majo1} we obtain that for any real number $M$ large enough,
$$
\PP^{\star}\left( \sqrt{n}\|\overline{\theta}_{n}({\cal K})-\theta^{\star}\|\geq M \right)\leq \PP^{\star}\left( W_{n}\geq \frac{\delta}{2CM({\cal K})}\right)+
\PP^{\star}\left( \sqrt{n}W_{n}\geq \frac{M\lambda}{2C}\right)
$$
where $M({\cal K})=\sup_{\theta\in{\cal K}} \|\theta\|$. This last equation together with Assumption $\mathbf{(A3)}$ gives the Theorem.
 \section{Proof of Theorem 4.1} \label{appb}
The proof follows the general methodology for model selection developed by  \citet{massart:2003}. 
To prove Theorem 4.1 and Theorem 4.2, we will use a concentration inequality 
we state now.
Let us introduce some notations. For any real function $f$,
denote
 $$
\G_{n} f =\frac{1}{\sqrt{n}}\sum_{i=1}^{n} \left[f\left(Y_{i}\right)-
\int f d\PP^{\star}\right].
$$
\begin{lemma}
\label{AdamBer}
Assume $\mathbf{(A4)}$. Let $\cal F$ be a class of real functions, and $F$ such that, for any $f\in{\cal F}$, $|f|\leq F$.
Assume that there exists $c(F)>0$ and $C(F)>0$ such that $\forall j=1,\ldots,k^{\star}$, $|g(j)| \leq C(F)$ where $g$ is defined by
$$
g\left(j\right)=\ln E_{\PP^{\star}}\left\{ \exp \left[2c(F)^{-1}|F(Y_{2})|\right]\vert S_{1}=j\right\}.
$$
Then there exist universal constants  $C_{1}$, $C_{2}$, $K_{1}$, $K_{2}$ and a constant $C^{\star}$ depending only on $Q^{\star}$ such that
\begin{multline*}
\PP^{\star} \left( \sqrt{n}\sup_{f\in {\cal F}} \G_{n} f  \geq K_{1} \sqrt{n} E_{\PP^{\star}}\left(\sup_{f\in {\cal F}} \G_{n} f\right) + C_{1}\tau \sqrt{nx} + C_{2}C^{\star}c(F)C(F)x \right) \\
\leq K_{2}\exp\left\{-x\right\}
\end{multline*}
where $\tau^{2}= \sup_{f\in{\cal F}}E_{\PP^{\star}}f^{2}\left(Y_{1}\right)$.
\end{lemma}

\noindent Proof of Lemma \ref{AdamBer}}

The lemma is an application of Theorem 7 in  \citet{adamczak:bednorz:12} to the stationary Markov chain $(X_{i})_{i\geq 1}=(S_{i},Y_{i})_{i\geq 1}$ and functions $f(s,y):=f(y)$.
 Then, with the notations of  \citet{adamczak:bednorz:12} we get that:
\begin{itemize}
\item
$m=1$,
\item
the small set $C$ is the whole space,
\item
the minorizing probability measure $\nu$ is that of $(\tilde{S}_{i},Y_{i})_{i\geq 1}$ with $(\tilde{S}_{i})_{i}$ i.i.d. with uniform distribution, and $\delta=\min_{i,j}Q^{\star}_{i,j}$.
\item
Since $C$ is the whole space, the return times $\sigma(i)=i$, so that $s_{i}(f)=f(Y_{i})$, thus the $\sigma^{2}$ of Theorem 7 is just
$\sup_{f} E_{\PP^{\star}}(f^{2}(Y_{1}))$.
\end{itemize}
Using the specific assumption of the lemma, 
 taking $\alpha=1$, we can apply  Corollary 1 of  \citet{adamczak:bednorz:12}, to get (with their notations again)
$$
a,b,c \leq C^{\star}c(F)C(F)
$$
for some constant $C^{\star}>0$ depending only on $\min_{i,j}Q^{\star}_{i,j}$.
$\Box$

For any $p\geq 2$, define
\begin{multline*}
S_{p}=\left\{\sum_{j=1}^{k^{\star}}\mu\left(j\right)f\left(\cdot-m_{j}\right), f\in {\cal F}_{p},\; |m_{j}|\leq M({\cal K}),\; \sum_{j=1}^{k^{\star}}\mu\left(j\right)=1,\;\right.\\
\left.\mu\left(j\right)\geq 0,\; j=1,\ldots,k^{\star}\right\},
\end{multline*}
so that $\widehat{s}_{p}\in S_{p}$. We now fix, for any $p\geq 2$, some $\tilde{s}_{p}\in S_{p}$ such that:
\begin{equation}
\label{norms}
\forall t\in S_{p},\;\|\sqrt{s^{\star}}-\sqrt{\tilde{s}_{p}}\|^{2}\leq 2 \|\sqrt{s^{\star}}-\sqrt{t}\|^{2}.
\end{equation}
For any $p\geq 2$ and any $\sigma >0$, define
$$
W_{p}\left(\sigma \right)=\sup_{t\in{\cal S}_{p}, \|\sqrt{t}-\sqrt{\tilde{s}_{p}}\|_{2}\leq \sigma} \G_{n}\left(\ln \left(\frac{s^{\star}+t}{s^{\star}+\tilde{s}_{p}}\right)\right),
$$
and let $L_{p}$ be an enveloppe function of $\{\ln \left(s^{\star}+t\right)-\ln\left(s^{\star}+\tilde{s}_{p}\right), t\in{\cal S}_{p}\}$.
Assume there exists functions $\psi_{p}$ such that  $\psi_{p}(x)/x$ is non increasing and for all $p\geq 2$ and $\sigma >0$,
\begin{equation}
\label{psip}
E_{\PP^{\star}}\left[W_{p}\left(\sigma \right)\right]\leq \psi_{p}\left(\sigma\right).
\end{equation}
Define $\sigma_{p}$ (depending also on $n$) as the unique solution of
\begin{equation}
\label{step4}
\psi_{p}\left(\sigma_{p}\right)=\sqrt{n} \sigma_{p}^{2}.
\end{equation}
Now we follow and adapt the proof of Theorem 7.11 in \citet{massart:2003}.
Let $p$ be such that $K(s^{\star},s_{p}) < +\infty$.
If $p'$ is such that $D(p') \leq D(p)$, then one gets, 
as in  \citet{massart:2003} p.241,
\begin{multline}
\label{step1}
K\left(s^{\star},\frac{s^{\star}+\widehat{s}_{p'}}{2}\right) 
-\frac{1}{2\sqrt{n}}\G_{n}\left(\log \left(\frac{s_{p}}{s^{\star}} \right)\right)
\\ \leq 
 K\left(s^{\star},s_{p}\right)+\text{pen}\left(p,n\right)-\frac{1}{\sqrt{n}}\G_{n}\left(\ln\left(\frac{s^{\star}+\widehat{s}_{p'}}{2s^{\star}}\right)\right)
-\text{pen}\left(p',n\right)+V_{p}
\end{multline}
where 
$$V_{p}=
\frac{1}{n}\sum_{i=1}^{n} \ln \left(\frac{\sum_{j=1}^{k^{\star}}\widehat{\mu}(j)f_{p}(Y_{i}-\widehat{m}_{j})}{\sum_{j=1}^{k^{\star}}\mu^{\star}(j)f_{p}(Y_{i}-m^{\star}_{j})} \right).
$$
Applying Lemma 4.23  in \citet{massart:2003} p. 139, for any positive $y_{p'}$:
$$
E^{\star}\left[\sup_{t\in{\cal S}_{p'}} \G_{n}\left(\frac{\ln \left(s^{\star}+t\right)-\ln\left(s^{\star}+\tilde{s}_{p'}\right)}{y_{p'}^{2}+\|\sqrt{t}-\sqrt{s_{p'}}\|_{2}^{2}}\right)\right]\leq 4\frac{\psi_{p'}\left(y_{p'}\right)}{y_{p'}^{2}}.
$$
Using  Lemma \ref{AdamBer}, the fact that $2y_{p'}\|\sqrt{t}-\sqrt{\tilde{s}_{p'}}\|_{2}\leq y_{p'}^{2}+\|\sqrt{t}-\sqrt{\tilde{s}_{p'}}\|_{2}^{2}$, and Lemma  7.26 p. 276  in \citet{massart:2003}, we obtain that
for some constant $C>0$, except on a set with probability less than $K_{2} \exp-(x_{p'}+x)$, for all $x>0$:
$$
\frac{1}{\sqrt{n}}\G_{n}\left(\frac{\ln \left(s^{\star}+\widehat{s}_{p'}\right)-\ln\left(s^{\star}+\tilde{s}_{p'}\right)
}{y_{p'}^{2}+\|\sqrt{\tilde{s}_{p'}}-\sqrt{\widehat{s}_{p'}}\|^{2}}\right)
\leq \frac{Cte}{y_{p'}} \left(\frac{\psi_{p'}\left(y_{p'}\right)}{y_{p'}\sqrt{n}}+\frac{T(L_{p'})(x_{p'}+x)}{ny_{p'}}+\sqrt{\frac{x_{p'}+x}{n}}
 \right).
$$
Here, $T(L_{p'})=C_{2}C^{\star}c(L_{p'})C(L_{p'})$.
Using again
 Lemma \ref{AdamBer} and Lemma  7.26 p. 276  in \citet{massart:2003}
we get that, for some constant $C>0$, except on a set with probability less than $K_{2} \exp-(x_{p'}+x)$, for all $x>0$:
$$
\frac{1}{\sqrt{n}}\frac{\G_{n}\left(\ln \left(s^{\star}+\tilde{s}_{p'}\right)-\ln\left(2s^{\star}\right)\right)}{y_{p'}^{2}+\|\sqrt{s^{\star}}-\sqrt{\tilde{s}_{p'}}\|^{2}}
\leq \frac{C}{y_{p'}} \left(\frac{T(L_{p'})(x_{p'}+x)}{ny_{p'}}+\sqrt{\frac{x_{p'}+x}{n}}
 \right).
$$
Now, using (\ref{norms}), we get
$$
\|\sqrt{\tilde{s}_{p'}}-\sqrt{\widehat{s}_{p'}}\|^{2}\leq \left[\|\sqrt{\tilde{s}_{p'}}-\sqrt{s^{\star}}\|+\|\sqrt{s^{\star}}-\sqrt{\widehat{s}_{p'}}\|\right]^{2}\leq 6\|\sqrt{s^{\star}}-\sqrt{\widehat{s}_{p'}}\|^{2}
$$
and
$$
\|\sqrt{s^{\star}}-\sqrt{\tilde{s}_{p'}}\|^{2}\leq 2\|\sqrt{s^{\star}}-\sqrt{\widehat{s}_{p'}}\|^{2}
$$
and we finally obtain that, for some other constant $C>0$ depending only on $\PP^{\star}$,
except on a set with probability less than $2K_{2} \exp-(x_{p'}+x)$, for all $x>0$:
$$
\frac{1}{\sqrt{n}}\G_{n}\left(\frac{\ln \left(s^{\star}+\widehat{s}_{p'}\right)-\ln\left(2s^{\star}\right)}{y_{p'}^{2}+\|\sqrt{s^{\star}}-\sqrt{\widehat{s}_{p'}}\|^{2}}\right)
\leq \frac{C}{y_{p'}} \left(\frac{\psi_{p'}\left(y_{p'}\right)}{y_{p'}\sqrt{n}}+\frac{T(L_{p'})(x_{p'}+x)}{ny_{p'}}+\sqrt{\frac{x_{p'}+x}{n}}
 \right).
$$
 Define for some constant $a$ to be chosen 
$$
y_{p'}=a^{-1}\sqrt{\sigma_{p'}^{2}+\frac{(x_{p'}+x)(1+T(L_{p'}))}{n}}.
$$
Then we can follow the proof of Theorem 7.11 in \citet{massart:2003} to obtain that, as soon as
\begin{equation}
\label{penmassart}
\text{pen}\left(p,n\right) \geq \kappa \left( \sigma_{p}^{2}+\frac{x_{p}(1+T(L_{p}))}{n}\right),
\end{equation}
one has for any $n\geq 2$, for some  real numbers $\kappa>0$ and $C>0$ depending only on $Q^{\star}$
$$
E_{\PP^{\star}}\left[h^{2}\left(s^{\star},\widehat{s}_{\widehat{p}}\right)\right] \leq 
C \left\{\inf_{p\geq 2}\left(K\left(s^{\star},s_{p}\right)+\text{pen}\left(p,n\right)+E_{\PP^{\star}}\left[V_{p}\right]
\right)+\frac{\Sigma}{n}
\right\}.
$$
But using the convexity of the Kullback-Leibler divergence to both arguments, we have, for any $p\geq 2$, $K\left(s^{\star},s_{p}\right)\leq K\left(f^{\star},f_{p}\right)$.
Thus to finish the proof of Theorem \ref{oracle}, one has to find functions $\psi_{p}$ verifying (\ref{psip}), evaluate $\sigma_{p}$ using (\ref{step4}), and evaluate
$T(L_{p})$.\\
Let us first prove that there exists constants $C,C'>0$ depending only on $\delta$ and $Q^{\star}$ such that, as soon as $\mathbf{(A4)}$  holds, for any $p\geq 2$,
\begin{equation}
\label{M(Lp)}
T\left(L_{p}\right) \leq C \ln \left(1+\frac{C'}{b_{p}^{\delta}}\right).
\end{equation}
First of all, we see that we can take
$$
L_{p}\left(y\right)=\ln \left( 1+\frac{1}{b_{p}\sqrt{2\pi}s^{\star}(y)} \right),
$$
with  $c(L_p)=2/\delta$, the function defined in Lemma \ref{AdamBer} is given by
$$
g\left(s\right)=\log \left[\sum_{j=1}^{k^{\star}} Q^{\star}_{s,j}\int \left( 1+\frac{1}{b_{p}\sqrt{2\pi}s^{\star}(u)}\right)^{\delta}f^{\star}(u-m^{\star}_{j})du\right]
$$
Under  $\mathbf{(A4)}$, on  gets that there exists constants $C>0$ depending only in $Q^{\star}$ and $\delta$ such that $g$ is bounded by the constant $ C \ln \left(1+\frac{C'}{b_{p}^\delta}\right)$ and (\ref{M(Lp)}) follows (for maybe another constant $C$).\\
%
To find functions $\psi_{p}$, we shall use \citet{doukhan:massart:rio:1995}. Since $(Y_{t})_{t\in\N}$ is geometrically ergodic,  Lemma 2 in \citet{doukhan:massart:rio:1995},  implies that, for some constant $C$ that depends only on $Q^{\star}$, for any real function $f$,
$$
\left\|  f  \right\|_{\beta}^{2} \leq C  \gamma(f) ( 1 + \log^+(\gamma(f)), \quad \gamma(f) = \int f^{2} ( 1 + \log^{+}|f| ) d\PP^{\star}  
$$
where $\left\| \cdot  \right\|_{\beta}$ is defined in \citet{doukhan:massart:rio:1995}. Now, since for all $x>0$, $x \ln^{+} x \leq x^{2}/e$,
$$
\left\|  f  \right\|_{\beta}^{2} \leq \frac{C}{e}\int |f|^{3}d\PP^{\star} \left(1+ \log^{+}(\frac{1}{e}\int |f|^{3}d\PP^{\star})\right). 
$$
Using Lemma 7.26 in \citet{massart:2003}, we thus get that for all $t\in{\cal S}_{p}$,
\begin{eqnarray*}
\left\|
\ln \left(s^{\star}+t\right)-\ln\left(s^{\star}+\tilde{s}_{p}\right)
\right\|_{\beta}^{2}
& \leq &\frac{1}{c^{2}}\|\sqrt{t}-\sqrt{\tilde{s}_{p}}\|_{2}^{2}
\end{eqnarray*}
for some constant $c>0$ that depends only on $Q^{\star}$, and the same trick leads to 
$$
H_{\beta}\left(u,\{\ln \left(s^{\star}+t\right)-\ln\left(s^{\star}+\tilde{s}_{p}\right), \;t\in{\cal S}_{p}\}\right)
\leq
H_{2}\left(c u,\{\sqrt{t},\;t\in{\cal S}_{p}\}\right)
$$
where $H_{\beta}\left(u,{\cal F}\right)$ is the bracketing entropy of a set ${\cal F}$ at level $u$ with respect to $\|\cdot \|_{\beta}$, that is the logarithm of the minimum of the number of brackets of $\|\cdot \|_{\beta}$-width $u$ needed to cover ${\cal F}$, and
$H_{2}\left(u,{\cal F}\right)$ is the bracketing entropy of a set ${\cal F}$ at level $u$ with respect to $\|\cdot \|_{2}$.\\
Let for any for $\sigma >0$ and $p\geq 2$
$$
\eta_{p} \left(\sigma\right)=\int_{0}^{\sigma/c}\sqrt{H_{2}\left(c u,\{\sqrt{t},\;t\in{\cal S}_{p}\}\right)}du.
$$
Using  Theorem 3 in \citet{doukhan:massart:rio:1995} we get
\begin{equation}
\label{maxineq}
E_{\PP^{\star}}\left[W_{p}\left(\sigma\right) \right] \leq A \eta_{p} \left(\sigma\right) \left[ 1 + \frac{\delta_{p}(1\wedge \epsilon (\sigma,n))}{\sigma} \right],
\end{equation}
where $ \epsilon (\sigma,n)$ is the unique solution of $x^{2}/B(x)=\eta_{p}^{2}(\sigma)/n\sigma^{2}$, $$B(x)=x+C(x-x\ln x)$$ for some constant $C$ that depends only on $Q^{\star}$, and $\delta_{p}$ is the function given by
$$
\delta_{p}\left(\epsilon\right)=\sup_{t\leq \epsilon}Q\left(t\right)\sqrt{B(t)}
$$
with for any $t$, $Q\left(t\right)\leq u$ iff $\PP^{\star}(H_{p}(Y_{1})> u)\leq t$. Here, $H_{p}$ is an envelope function of
$\{\ln \left(s^{\star}+t\right)-\ln\left(s^{\star}+\tilde{s}_{p}\right), \;\|\sqrt{t}-\sqrt{\tilde{s}_{p}}\|\leq \sigma,\;t\in{\cal S}_{p}\}$. Taking $H_{p}=L_{p}$ one gets easily
$$
Q\left(t\right)\leq \ln \left(1+\frac{1}{tb_{p}\sqrt{2\pi}}\right),
$$
so that
$
\delta_{p}\left(\epsilon\right)\leq \sup_{t\leq \epsilon}h_{p}(t)
$
with 
$$
h_{p}\left(t\right)=\ln \left(1+\frac{1}{tb_{p}\sqrt{2\pi}}\right)\sqrt{t+C(t-t\ln t)}.
$$
The variations of $h_{p}$ imply that there exists a universal constant $b$ such that as soon as $b_{p}\leq b$, $h(t)$ is increasing on $(0,1)$, so that
$$
\delta_{p}\left(\epsilon \wedge 1\right)=h_{p}\left(\epsilon \wedge 1\right) \leq \tilde{h}_{p}\left(\epsilon \wedge 1\right) 
$$
with
$$
 \tilde{h}_{p}\left(t\right) =C \ln \left(\frac{1}{b_{p}} \right)\sqrt{t} |\ln t | \left(\sqrt{|\ln t |}\wedge 1\right),
$$
 for some universal constant $C$. Using \citet{MM08}, we get that for some fixed constant $K$, for all $u>0$,
$$
H_{2}\left(u,\{\sqrt{t},\;t\in{\cal S}_{p}\}\right)\leq k^{\star}p\left[3 \ln \left(\frac{1}{u\wedge 1}\right) + \frac{3}{4} \ln \left(\frac{1}{b_{p}}\right)
+ \ln A_{p} + K\right]+ \ln \left(k^{\star}p\right).
$$
Using the fact that for all $\epsilon\in]0,1]$,
$$
\int_{0}^{\epsilon}\sqrt{\ln \left(\frac{1}{x}\right)}dx \leq \epsilon \left\{\sqrt{\ln \left(\frac{1}{\epsilon}\right)}+\sqrt{\pi}\right\}
$$
and since
$$
\eta_{p} \left(\sigma\right)
=
\frac{1}{c}\int_{0}^{\sigma}\sqrt{H_{2}\left(c^{2} u,\{\sqrt{t},\;t\in{\cal S}_{p}\}\right)}du
$$
we get that for some other fixed constant $K$ and all $\sigma>0$
\begin{equation}
\label{phip}
\eta_{p} \left(\sigma\right)
\leq
\frac{\sigma}{c}t_{p}\left(\sigma\right)
\end{equation}
with
$$
t_{p}\left(\sigma\right)=\sqrt{k^{\star}p}\left[3\sqrt{\ln \left(\frac{1}{\sigma\wedge 1}\right)}+\sqrt{\frac{3}{4} \ln \left(\frac{1}{b_{p}}\right)
+ \ln A_{p} + K}\right]+\sqrt{\ln (k^{\star}p)}.
$$
Now, one may use the upper bound \eqref{phip} to upper bound $\epsilon (\sigma , n)$, and we get that for some universal constant $C$, 
$$
\epsilon \left(\sigma , n\right)\leq C \frac{t_{p}^{2}\left(\sigma\right)}{n}\ln \left[  \frac{t_{p}^{2}\left(\sigma\right)}{2n}\right].
$$
Then we may set
$$
\psi_{p}\left(\sigma\right)=\frac{\sigma}{c}t_{p}\left(\sigma\right)\left[ 1 + \frac{\tilde{h}\left(C \frac{t_{p}^{2}\left(\sigma\right)}{n}\ln \left[  \frac{t_{p}^{2}\left(\sigma\right)}{2n}\right]\wedge 1\right)}{\sigma}\right].
$$
(\ref{psip}) holds, $\psi_{p}(x)/x$ is indeed non increasing, and if $\sigma_{p}$ is the unique solution of (\ref{step4}), we obtain that for some constant $C$ depending only on $\PP^{\star}$,
as soon as $b_{p}\leq b$,
\begin{equation}
\label{sigmap}
\sigma_{p}^{2}\geq \frac{C}{n}k^{\star}p \left[\ln n + \ln\left(\frac{1}{b_{p}}\right)+\ln A_{p}\right].
\end{equation}
\section{ Proof of Theorem 4.2} \label{subsec:Vp}


For simplicity's sake we denote in the following $ \mathcal H_{loc} (\beta):= \mathcal H_{loc}(\beta, \gamma,{\cal P})$.
Set  $p = p_0\lfloor (n/\log n)^{1/(2\beta + 1)} (\log n)^{4 \beta / (2\beta +1)} \rfloor$ with $p_0>0$ fixed which we shall determine later,
$b_p = b_0 (\log p)^{2}/p$ for some positive $b_0$ and $A_{p}=a_{0}|\log b_{p}|$ for some positive $a_{0}$. 
The approximating  $f_{p}\in {\cal F}_{p}$ 
$$ f_p(y) = \sum_{i=1}^p \pi_i \varphi_{b_p}( y-\alpha_i)
$$
is taken from \citet{kruijer:rousseau:vdv:09}. Let $\ell_{j}^{\star}$ denote the $j$-th derivative of $\log f^{\star}$.
A simple modification in the proof of Lemma 4 of \citet{kruijer:rousseau:vdv:09} gives that for any $H$ and any $\tilde H$ with $H > \tilde H + 3 \beta$, there exists
$\tilde B$ such that if 
\begin{multline*}
D_{p}:= \left\{ y\;:\; f^{\star}(y-m) \geq b_p^{\tilde{H}}, |\ell_{j}^{\star}(y-m) | \leq \tilde{B} b_p^{-j} |\log p|^{-j/2}, \,j\leq \beta ,\right.\\
\left. |L(y-m) |   \leq \tilde{B} b_p^{-\beta} |\log p|^{-\beta/2}, \, \forall 0 \leq m \leq 2m^{\star}_k\right\}
\end{multline*}
then, for all $y \in D_{p}$ and all $0 \leq m \leq m^{\star}_k$
 \begin{equation} \label{lem:KRVDV:1}
f_{p}(y-m) =  f^{\star}(y-m) ( 1 + O(R(y-m)b_p^\beta)) + O((1+ R(y-m))b_p^{H-\tilde{H}}),
 \end{equation}
where the function $R(y)$ is a linear combination of $L(y) $ and of the functions $|\ell_{j}^{\star}(y)|^{\beta /j}$, $j \leq \beta$, 
and where the constants entering the terms $O(.)$ depend on $\mathcal
H_{loc}(\beta)$, $\tilde B$, $H$ and $\tilde H$. Note that since the functions $l_j^{\star}$ are bounded by polynomials, there exists a constant $C$ such that 
 $ |R(y-m)| \leq C(1+ R(y))$, $\forall 0 \leq m \leq 2m^*_k$. 
In the following we fix $\tilde H > 4 \beta + 2 \gamma$ and $H >\tilde H +
3\beta $.
Moreover, Lemma 4 in \citet{kruijer:rousseau:vdv:09} implies
 \begin{equation} 
 \label{equakul}
 K( f^{\star} , f_p) \lesssim b_p^{2\beta}, \quad \int f^{\star} \left( \log f^{\star} - \log f_p \right)^2 (y) dy \lesssim b_p^{2\beta}.
 \end{equation}
 Here and further, $\lesssim$ will denote an upper bound up to a constant, where the constant entering the upper bound depends only on  $\mathcal H_{loc}(\beta)$.
 Throughout the proof $C$ denotes a generic constant depending only on $H_{loc}(\beta)$ and $Q^{\star}$. 
 
 First of all, with such choices of $p$, $b_{p}$, and $A_{p}$, using  Theorem \ref{oracle} and \eqref{equakul}, there remains to prove that   $E_{\PP^{\star}}[V_{p}]\lesssim b_p^{2\beta}$ or equivalently
 $v_{n}E_{\PP^{\star}}[V_{p}]\lesssim 1$ with  $v_{n}= n^{\frac{2\beta }{2\beta+1}}(\log n)^{-6\beta /(2\beta + 1) }$.
For any $\theta$ and any $y$, set
  $$
  w_{p}\left(\theta,y\right)=  \log\left(  \frac{  \sum_{j=1}^{k^{\star} } \mu(j) f_p( y -  m_j) }{ \sum_{j=1}^{k^{\star}} \mu^{\star}(j) f_p( y-  m_{j}^{\star}) }\right).
   $$
First note that 
 \begin{equation} \label{ineq:LLR}
  \log \left( \frac{ \sum_{j=1}^k  \mu_j^{\star} f_p( y - m_j) }{ \sum_{j=1}^k \mu_j^{\star} f_p( y -  m_j^{\star}) } \right)
   \leq \max_j \left\{ \frac{( | y-m_j^{\star}|  +A_{p}) |m_j - m_j^{\star}|  }{   b_p^2 }\right \}.
  \end{equation}
Thus we can bound
 \begin{equation*}
 \begin{split}
 &   \frac{v_n}{n} E_{\PP^{\star}}\left[\sum_{i=1}^{n} w_{p}(\widehat{\theta},Y_{i})\ind_{\| \theta - \widehat \theta\| > M_0 \sqrt{\log n/n} } \right] \\
 &\leq   \frac{v_n}{n}\max_j   \sum_{i=1}^n  E_{\PP^{\star}}\left[\ind_{\| \theta - \widehat \theta\| > M_0 \sqrt{\log n/n} }\left( \frac{( | Y_i-m_j^{\star}|  +A_{p}) |\hat m_j - m_j^{\star}|  }{   b_p^2 }  + \frac{ | \hat \mu_j - \mu_j^{\star} | }{ \mu_j^{\star} } \right) \right] \\
 &\lesssim \left(  \frac{ v_{n}\log p  }{ b_p^2 \sqrt{n} } + \frac{v_n}{\sqrt{n}} \right) \PP^\star \left[ \| \theta - \widehat \theta\| > M_0 \sqrt{\log n/n} \right] \\
 &= o(1)
  \end{split} 
  \end{equation*}
  by Theorem \ref{theotailstheta} and choosing $M_0 =1/\sqrt{c^{\star}}$. 
  \\
  Set now $H_1 > 3
+ 2\beta$,  $C_{p,1}=D_p^c \cap \{ |y| \leq H_1 \log (1/b_p)\tau^{-1} \}$ and  $C_{p,2}=D_p^c \cap \{ |y| > H_1 \log (1/b_p)\tau^{-1} \}$.
  Using  \eqref{ineq:LLR} we get, for all $i=1, \cdots, n$,
\begin{equation*}
\begin{split}
 E_{\PP^{\star}}\left[\ind_{C_{p,1}}(Y_i) w_{p}(\widehat{\theta},Y_{i})\ind_{\| \theta - \widehat \theta\| \leq M_0 \sqrt{\log n/n}} \right]
 &\lesssim 
\frac{   (\log p)^{3/2} }{   b_p^2\sqrt{n} } \int_{C_{p,1}} s^{\star}(y) dy \\
 &\leq \frac{ (\log p)^{3/2}   }{  b_p^2 \sqrt{n}  } b_p^{2\beta+ \gamma}  \lesssim v_n^{-1}
\end{split}
\end{equation*}
as soon as $\gamma > (3/2 - \beta)_+$, where the last inequality comes from an adaptation of Lemma 2 in \citet{kruijer:rousseau:vdv:09}, using the moment conditions \eqref{moment}. We also have
\begin{equation} \label{Epc}
\begin{split}
  E_{\PP^{\star}}\left[\ind_{C_{p,2}}(Y_i) w_{p}(\widehat{\theta},Y_{i})\ind_{\| \theta - \widehat \theta\| \leq M_0 \sqrt{\log n/n}}\right]
 &\lesssim
\frac{ (\log p)^{3/2}   }{  b_p^2 \sqrt{n}  } \int_{C_{p,2}}|y| s^{\star}(y)dy \\
 &\lesssim  \frac{ (\log p)^{3/2}b_p^{H_1/2} }{  b_p^2 \sqrt{n}  } \lesssim v_n^{-1}
\end{split}
\end{equation}
since $H_1 > 3
+ 2\beta$, where the last inequality comes from the tail condition \eqref{cond:tail}. There thus remains to prove that
$$
v_{n}E_{\PP^{\star}}\left[\frac{1}{n}\sum_{i=1}^{n} w_{p}(\widehat{\theta}_{n},Y_{i})\ind_{D_{p}}(Y_i)\ind_{\| \theta - \widehat \theta\| \leq M_0 \sqrt{\log n/n}}\right] \lesssim 1.
$$
    We shall use:
 \begin{multline}
\label{majointegr}
v_{n}E_{\PP^{\star}}\left[\frac{1}{n}\sum_{i=1}^{n} w_{p}(\widehat{\theta}_{n},Y_{i})\ind_{D_{p}}(Y_i)\ind_{\| \theta - \widehat \theta\| \leq M_0 \sqrt{\log n/n}}\right]\\
\leq\int_{0}^{+\infty}\PP^{\star}\left(\frac{v_{n}}{n}\sum_{i=1}^{n} w_{p}(\widehat{\theta}_{n},Y_{i})\ind_{D_{p}}(Y_i)\ind_{\| \theta - \widehat \theta\| \leq M_0 \sqrt{\log n/n}}\geq x\right)dx.
\end{multline}   
  Notice now that
  \begin{multline*}
 \PP^{\star}\left(\frac{v_{n}}{n}\sum_{i=1}^{n} w_{p}(\widehat{\theta}_{n},Y_{i})\ind_{D_{p}}(Y_i)\ind_{\| \theta - \widehat \theta\| \leq M_0 \sqrt{\log n/n}}\geq x\right)\leq \\
  \PP^{\star}\left(\sqrt{n}\|\widehat{\theta}_{n}-\theta^{\star}\|\geq M(x)\right) +       \PP^{\star}\left(\sup_{\sqrt{n}\|\theta-\theta^{\star}\|\leq M(x) \wedge M_0 \sqrt{\log n/n}}\frac{v_{n}}{\sqrt{n}}\G_{n}( w_{p}(\theta,\cdot)\ind_{D_{p}}(\cdot))\geq \frac{x}{2}\right)
          \end{multline*} 
as soon as
  \begin{equation}
 \label{majoesp1}
 v_{n}\sup_{\sqrt{n}\|\theta-\theta^{\star}\|\leq M(x)\wedge M_0 \sqrt{\log n}}E_{\PP^{\star}}[w_{p}(\theta,Y_{1})\ind_{D_{p}}(Y_1)]\leq \frac{x}{2}.
 \end{equation}
 If moreover
 \begin{equation}
 \label{majoesp2}
\frac{v_{n}K_{1}}{\sqrt{n}}E_{\PP^{\star}}\left(\sup_{\sqrt{n}\|\theta-\theta^{\star}\|\leq M(x)\wedge M_0 \sqrt{\log n}}\G_{n}( w_{p}(\theta,\cdot)\ind_{D_{p}}(\cdot))\right)\leq \frac{x}{4},
 \end{equation}
where $K_1$ is defined in Lemma \ref{AdamBer}, Appendix \ref{appb},  using Theorem \ref{theotailstheta} and Lemma \ref{AdamBer} we get, for large enough $x$, with $M(x) = x^{1/4}$,
  \begin{multline*}
 \PP^{\star}\left(|\frac{v_{n}}{n}\sum_{i=1}^{n} w_{p}(\widehat{\theta}_{n},Y_{i})\ind_{D_{p}}(Y_i)\ind_{\| \theta - \widehat \theta\| \leq M_0 \sqrt{\log n/n}}|\geq x\right)\leq \\
2\exp\left( -\frac{nx}{v_{n}C_{n}(x)}\right)+ 2\exp\left(-\frac{nx^{2}}{v_{n}^{2}\tau_{n}(x)^{2}}\right) +
8\exp\left(-c^{\star}x^{1/2}\right) 
\end{multline*}
with
$$
\tau_{n}(x)^{2}=16C_{1}^{2}\sup_{\sqrt{n}\|\theta-\theta^{\star}\|\leq M(x) \wedge M_0 \sqrt{\log n}}E_{\PP^{\star}}[w_{p}^{2}(\theta,Y_{1})\ind_{D_{p}}(Y_{1})],
$$
$$
C_{n}(x)=4C_{2}C^{\star}c\left(W_{n,p,x}\right)C\left(W_{n,p,x}\right),
$$
where $W_{n,p,x}$ is such that
$$
\sup_{\sqrt{n}\|\theta-\theta^{\star}\|\leq M(x) \wedge M_0 \sqrt{\log n}} w_{p}(\theta,\cdot) \ind_{D_p} \leq W_{n,p,x}(\cdot).
$$
For instance we may take
 $$
 W_{n,p,x}\left(y\right)=\log p + C\frac{(|y|+A_{p})M(x)}{\sqrt{n}b_{p}^{2}},
$$
leading, by choosing $c( W_{n,p,x})=C\frac{M(x)}{\sqrt{n}b_{p}^{2}\log n}$ , to
\begin{equation}
\label{cn}
C_{n}(x)=C\frac{A_{p}x^{1/2}}{\sqrt{n}b_{p}^{2}\log n}.
\end{equation}
For any $\theta$
  set
$$
s_{p,\theta}(y)=\sum_{j=1}^{k^{\star}}\mu (j) f_{p}(y-m_{j})\;\;\text{and}\;\;s^{\star}_{\theta}(y)=\sum_{j=1}^{k^{\star}}\mu (j) f^{\star}(y-m_{j}).
$$
We consider the following decomposition,
\begin{equation} \label{decomp:wtheta}
 \log \left( \frac{ s_{p,\theta}( y ) }{ s_{p, \theta^{\star}}(y) } \right) =  \log \left( \frac{ s_{p,\theta}( y ) }{ s_{ \theta}^{\star} (y) } \right)+  \log \left( \frac{  s_{ \theta}^{\star} (y) }{ s^{\star}(y) } \right)+ \log \left( \frac{ s^{\star}(y) }{ s_{p, \theta^{\star}}(y) } \right).
\end{equation}
 The first and third terms of \eqref{decomp:wtheta} are treated similarly. \eqref{lem:KRVDV:1} gives  that  for any $\theta$, over $D_p$,
\begin{equation}
\label{term1et3}
\left| \log \left(\frac{ s_{p,\theta}( y ) }{ s_{ \theta}^{\star} (y) }\right) \right| \lesssim R(y) b_p^\beta.
\end{equation}
For the second term, 
since $f^{\star}\in H_{loc}(\beta, \gamma, {\cal P})$ with $\beta \geq 1/2$,
$$|\log f^{\star}(y - \widehat m_j)- \log f^{\star}(y - m_j^{\star}) | \leq L(y-m_j^{\star}) |\widehat m_j - m_j^{\star}|^{\beta \wedge 1}.$$
Moreover,  if $y\in D_{p}$, and $\sqrt{n}\|\theta-\theta^{\star}\|\leq M(x) \wedge M_0 \sqrt{\log n}$, then for large enough $n$,
$$
L(y-m_j^{\star}) |\widehat m_j - m_j^{\star}|^{\beta \wedge 1}\leq 1
$$
so that
we have, for $\theta$ such that $\sqrt{n}\|\theta-\theta^{\star}\|\leq M(x) \wedge M_0 \sqrt{\log n}$, over $D_{p}$, for large enough $n$,
 \begin{multline}
 \label{term2}
 \left| \log \left( \frac{ s_{\theta}^{\star}( y ) }{ s^{\star}(y) } \right) \right| 
 \lesssim  \sum_j \frac{ |\mu_j - \mu_j^{\star}|}{ \mu_j^{\star}} 
+ 
 \sum_j |m_j - m_j^{\star}|^{\beta \wedge 1}  L(y-m_j^{\star})\\
 \lesssim \frac{ M(x) }{ \sqrt{n} }+ (n^{-1/2}M(x))^{\beta \wedge 1}  \sum_j  L(y-m_j^{\star}).
\end{multline}
Thus, using the fact that $\beta \geq 1/2$, for large enough $x$, 
\begin{equation}
\label{taun}
\sup_{\sqrt{n}\|\theta-\theta^{\star}\|\leq M(x) \wedge M_0 \sqrt{\log n}}E_{\PP^{\star}}[w_{p}^{2}(\theta,Y_{1})\ind_{D_{p}}(Y_{1})]= O(M(x)^{2}b_p^{2\beta}).
\end{equation}
 \eqref{cn} and \eqref{taun} give that, for all $\beta \geq 1/2$, for large enough $x$,
  $$ \frac{nx}{v_nC_n(x) } \gtrsim  x^{1/2}n^{3/2}b_p^{2\beta+2} \gtrsim x^{1/2}(\log n)^{3(2\beta + 2)/(2\beta +1)},
  \;\;  
  \frac{ nx^2}{ v_n^2 \tau_n^2 (x) } \gtrsim x^{1/2} n ,$$
  so that for large enough $x$, 
 \begin{equation}
 \label{majoprob}
 \PP^{\star}\left(\frac{v_{n}}{n}\sum_{i=1}^{n} w_{p}(\widehat{\theta}_{n},Y_{i})\ind_{D_{p}}(Y_i)\ind_{\| \theta - \widehat \theta\| \leq M_0 \sqrt{\log n/n}}\geq x\right)\lesssim \exp\left(-C x^{1/2}\right) \end{equation}
as soon as  \eqref{majoesp1} and  \eqref{majoesp2} hold for large enough $x$.

We now prove \eqref{majoesp1}. 
\begin{equation*} 
\begin{split}
E_{\PP^{\star}}[w_{p}(\theta,Y_{1})\ind_{D_{p}}(Y_{1})] &= \int_{D_p} (s^{\star}(y)  - s_{p,\theta^{\star}}(y) ) \log \left( \frac{ s_{p,\theta}( y ) }{ s_{p, \theta^{\star}}(y)} \right)dy - K(s_{p,\theta^{\star}}, s_{p,\theta}) \\
&\leq \int_{D_p} (s^{\star}(y)  - s_{p,\theta^{\star}}(y) ) \log \left( \frac{ s_{p,\theta}( y ) }{ s_{p, \theta^{\star}}(y) } \right)dy.
\end{split}
\end{equation*}
 Moreover, \eqref{lem:KRVDV:1} and \eqref{term1et3} give  that
\begin{multline*}
 \int_{D_p} |s^*(y)  - s_{p,\theta^{\star}}(y) | \left| \log \left(\frac{ s_{p,\theta}( y ) }{ s_{ \theta}^{\star} (y) }\right) \right|dy\\
  \lesssim b_p^{\beta}h(s^{\star}, s_{p, \theta^{\star}}) \left(\int_{D_p} |R(y) |^2 (s^{\star}(y) + s_{p,\theta^{\star}}(y))dy \right)^{1/2} 
 \lesssim b_p^{2\beta}
 \end{multline*}
using \eqref{equakul}.
Also, \eqref{lem:KRVDV:1} and \eqref{term2} give  that for $\theta$ such that $\sqrt{n}\|\theta-\theta^{\star}\|\leq M(x) \wedge M_0 \sqrt{\log n}$, 
 \begin{equation*}
  \int_{D_p} |s^{\star}(y)  - s_{p,\theta^{\star}}(y) | \left| \log \left( \frac{ s_{\theta}^{\star}( y ) }{ s^{\star}(y) } \right) \right| dy 
 \lesssim b_p^{2\beta}M(x)^{\beta \wedge 1}
\end{equation*}
so that for $\beta \geq 1/2$, uniformly for  $\theta$ such that $\sqrt{n}\|\theta-\theta^{\star}\|\leq M(x)\wedge M_0 \sqrt{\log n}$, 
$$E_{\PP^{\star}}[w_{p}(\theta,Y_{1})\ind_{D_p}(Y_1)] = O(M(x)b_p^{2\beta})$$
and \eqref{majoesp1} holds for large enough $x$.

  To Prove  \eqref{majoesp2} we use \eqref{decomp:wtheta}. We first control 
   $$ \mathbb E_{\PP^{\star}} \left( \sup_{\sqrt{n}\|\theta - \theta^{\star}\|\leq M(x)\wedge M_0\sqrt{\log n}}  \mathbf G_n (\ind_{D_p}  \log ( s_{p, \theta}/s^{\star}_{\theta}))\right).$$ 
   Using  \eqref{term1et3}, we can bound on $D_{p}$,
   $$\left| \log \left( \frac{s_{p, \theta}}{s^{\star}_{\theta}}(y) \right)\right|  \lesssim  |R(y)| b_p^{\beta} \lesssim (\log p)^{-1/2}\leq 1$$
for $n$ large enough, uniformly over $\sqrt{n}\|\theta - \theta^{\star}\|\leq M(x)\wedge M_0\sqrt{\log n}$. 
Also, $\|f\|_{2,\beta}^2 \lesssim \int f^2 ( 1 + \log^+|f|)(y) dy\lesssim  \|f\|_2^2 $, for any $f$ in the form 
    $\log \left(  s_{p, \theta}/s^{\star}_{\theta}\right)$. We denote
    $$\varphi_1(\sigma) = \int_{0}^\sigma \sqrt{\mathcal H ( u , S_{n,p,1}(\sigma) , \| . \|_{2} ) } du $$
 with 
 $$ S_{n,p,1}(\sigma,x) = \{ \log \left(  s_{p, \theta}/s^{\star}_{\theta}\right), \sqrt{n}\|\theta - \theta^{\star}\|\leq M(x)\wedge M_0\sqrt{\log n},  \| \log (  s_{p, \theta}/s^{\star}_{\theta}) \|_2 \leq \sigma\}.$$
  Then 
for all $y \in D_p$, since $|y| \leq  A_{p}$, and for all $|m_j- m_j'| \leq \eta$,
  \begin{equation*}
  \begin{split}
  f_{p}(y - m_j') &=\sum_{l=1}^p \pi_l \varphi_{b_p}( y - m_j - \alpha_l) e^{  - \frac{ (m_j - m_j^{'})^2 }{ 2  b_p^2} + \frac{ ( y - m_j - \alpha_l)( m_j - m_j^{'})}{ b_p^2 }} \\
  & \leq f_p(y-m_j)  e^{   \frac{ (|y| + m_{k^{\star}} +A_p )\eta }{ b_p^2 }}  \\ 
  &\leq  f_p(y-m_j)  e^{   \frac{ 3A_p \eta }{ b_p^2 }}:= f_U(y-m_j)\\
  &\geq f_p(y-m_j)  e^{   -\frac{ \eta^2 }{ 2  b_p^2} - \frac{ 3A_p \eta }{ b_p^2 }}  := f_L(y-m_j)
  \end{split}
  \end{equation*}
  and 
   \begin{equation} \label{bound:f*}
  \begin{split}
  f^{\star}(y - m_j') &\leq f^{\star}( y - m_j) e^{ \eta^{\beta\wedge 1} \sup_{|m-m_j|<\eta} |\tilde \ell(y-m)|}\\
  & \geq f^{\star}( y - m_j) e^{ -\eta^{\beta\wedge 1}  \sup_{|m-m_j|<\eta} |\tilde \ell(y-m)|}
  \end{split}
  \end{equation}
where $\tilde \ell (y-m ) = \ell_1( y-m) $ if $\beta > 1$ and $\tilde \ell(y-m) = L(y-m_j)$ if $\beta \leq 1$, so that a bracket for $ \log ( s_{p,\theta'} / s^*{\theta'} )\ind_{D_{p}} $ is given on $D_{p}$ by 
\begin{equation*}
\begin{split}
U_{p,\theta} &:=   \left( \frac{ 3A_p \eta}{ b_p^2 } +  \eta^{\beta\wedge 1} \sup_{|m-m_j|<\eta} |\tilde \ell(y-m)|\right)+ \log (1 + \eta\sum_{ j =1}^{k^{\star}} \mu(j)^{-1} ) \\
 L_{p,\theta} &:= -\frac{3A_p \eta}{ b_p^2 } -  \eta^{\beta\wedge 1} \sup_{|m-m_j|<\eta} |\tilde \ell(y-m)|+ \log (1 - \eta\sum_{ j =1}^{k^{\star}} \mu(j)^{-1} ),
 \end{split}
 \end{equation*}
Thus if $u>0$ and $\eta \leq \eta_0 (u^{1/2} b_p^2 / A_p \wedge u^{(\beta\wedge 1)/2})$ with $\eta_0 >0$ small enough,
$$\int_{D_{p}} (U_{p,\theta} -L_{p,\theta})^2(y) s^{\star}(y) dy \leq u,$$
 so that 
 \begin{equation*}
 \varphi_1(\sigma) \lesssim \sigma \sqrt{ \log^+( 1/ \sigma ) + \log (nM(x)) }.
 \end{equation*}
 Moreover for all $\|\theta - \theta^{\star}\|\leq M_0\sqrt{\log n} /\sqrt{n}$, \eqref{lem:KRVDV:1} implies that 
 $$\| \log ( s_{p,\theta'} / s^{\star}_{\theta'} ) \|_2^2 \leq b_p^{2\beta } C.$$
 Therefore 
 using Theorem 2 of \citet{doukhan:massart:rio:1995} and the fact that the chain is geometrically ergodic,
 we obtain that 
 \begin{equation}\label{EGn1}
 \begin{split}
 \mathbb E_{\PP^{\star}} \left( \sup_{\sqrt{n}\|\theta - \theta^{\star}\|\leq M(x)\wedge M_0\sqrt{\log n}}  \mathbf G_n (\ind_{D_p}  \log ( s_{p, \theta}/s^{\star}_{\theta}))\right) 
 &\lesssim  b_p^\beta (\log n + \log M(x))\\
& \leq x/8,
 \end{split}
 \end{equation}
 for $x\geq 1$ and large enough $n$.
We now study 
$$ \mathbb E_{\PP^{\star}} \left( \sup_{\sqrt{n}|\theta - \theta^{\star}\|\leq M(x)\wedge M_0\sqrt{\log n}}  \mathbf G_n \left[\ind_{D_p}  \log \left( \frac{  s_{ \theta}^{\star} (y) }{ s^{\star}(y) } \right)\right]\right).$$
Using \eqref{bound:f*}, if $\sqrt{n}|\theta - \theta^{\star}\|\leq M(x)\wedge M_0\sqrt{\log n}$, 
\begin{equation*}
\left|  \log \left( \frac{  s_{ \theta}^{\star} (y) }{ s^{\star}(y) } \right) \right| \lesssim b_p^{-(\beta \wedge 1)}\sqrt{\log n}/\sqrt{n} = o(1), 
\end{equation*}
so that 
\begin{equation*}
\begin{split}
\left\| \log \left( \frac{  s_{ \theta}^{\star} (y) }{ s^{\star}(y) } \right) \right\|_{2,\beta}^2 &\lesssim \left\| \log \left( \frac{  s_{ \theta}^{\star} (y) }{ s^{\star}(y) } \right) \right\|_2^2 \\
&\lesssim  \max_j ( \mu_j/\mu_j^{\star} -1)^2 + \max_j \int s^{\star}(y) ( \log f^{\star}(y-m_j) - \log f^{\star}(y-m_j^{\star}))^2dy \\
&\lesssim (M(x)^2/n)^{\beta \wedge 1}.
\end{split}
\end{equation*}
Hence using the same tricks as before and applying Theorem 2 of \citet{doukhan:massart:rio:1995} we obtain that   for large enough $n$,  
 \begin{equation}
 \mathbb E_{\PP^{\star}} \left( \sup_{\sqrt{n}\|\theta - \theta^{\star}\|\leq M(x)\wedge M_{0}\sqrt{\log n}}  \mathbf G_n (\ind_{D_p}  \log ( s_{p, \theta}/s^{\star}_{\theta}))\right)
 \lesssim M(x) n^{-(\beta \wedge1 )/2}\sqrt{\log n} = o(x \sqrt{n}/v_n)
 \end{equation}
    for all $x$ and \eqref{majoesp2} is satisfied.
 \\
  Finally, \eqref{majoprob} holds, which, together with \eqref{majointegr}
ends the proof of Theorem \ref{sadaptive}.

\bibliographystyle{apalike}
\bibliography{biblioHMM.bib}

\end{document}